\documentclass[11pt,twoside]{article}
\usepackage{amssymb}
\usepackage{amsmath}
\usepackage{latexsym}
\usepackage{amsthm}
\usepackage{epsfig}
\usepackage{mathptmx}
\theoremstyle{plain}

\def\Proof{\ifdim\lastskip<\aline\removelastskip\skipaline\fi
\noindent\it Proof. \rm}

\makeatletter
\def\thmhead@plain#1#2#3{%
  \thmname{#1}\thmnumber{\@ifnotempty{#1}{ }#2}%
  \thmnote{ \the\thm@notefont(#3)}}
\let\thmhead\thmhead@plain
\def\swappedhead#1#2#3{%
  \thmnumber{#2}\thmname{\@ifnotempty{#2}{. }#1}%
  \thmnote{ \the\thm@notefont(#3)}}
\makeatother

\theoremstyle{definition} 
 \newtheorem{definition}{Definition}[section]

\theoremstyle{plain}      
 
 \newtheorem{theorem}[definition]{Theorem}

\newcommand{\coco}{\vspace{-1.4mm}}
\begin{document}

\title{Deformation quantization: genesis, developments \\
and metamorphoses\thanks{This is an extended version of the 
presentation by D.S. on May 30, 2001, at the RCP in Strasbourg, 
completed November 27, 2001. Published in \textsl{Deformation quantization}
(G. Halbout, ed.), IRMA Lectures in Math. Theoret. Phys. \textbf{1}, 
9--54, Walter de Gruyter, Berlin 2002.}}

\author{Giuseppe Dito and Daniel Sternheimer\\\vspace*{-1mm}
{\small {\it Laboratoire Gevrey de Math\'ematique Physique, CNRS \ UMR 5029}}\\
\vspace*{-1mm}
{\small{\it D\'epartement de Math\'ematiques, Universit\'e de Bourgogne}}\\
\vspace*{-1mm}
{\small{\it BP 47870, F-21078 Dijon Cedex, France.}}\\
{\small \texttt{Giuseppe.Dito@u-bourgogne.fr, 
Daniel.Sternheimer@u-bourgogne.fr}}}
\markboth{\textsc{Giuseppe Dito and Daniel Sternheimer}}
{\textsc{Deformation Quantization: Genesis, Developments and Metamorphoses}}
\date{\small November 27, 2001}
\maketitle

\vspace*{-10mm}

\begin{abstract}
We start with a short exposition of developments in physics and
mathematics that preceded, formed the basis for, or accompanied, the
birth of deformation quantization in the seventies. We indicate how the
latter is at least a viable alternative, autonomous and conceptually more 
satisfactory, to conventional quantum mechanics and mention related questions,
including covariance and star representations of Lie groups. We sketch
Fedosov's geometric presentation, based on ideas coming from index theorems,
which provided a beautiful frame for developing existence and classification 
of star-products on symplectic manifolds. We present Kontsevich's formality, 
a major metamorphosis of deformation quantization, which implies existence
and classification of star-products on general Poisson manifolds and has
numerous ramifications. Its alternate proof using operads gave a new
metamorphosis which in particular showed that the proper context is that of 
deformations of algebras over operads, while still another is provided by
the extension from differential to algebraic geometry. In this panorama
some important aspects are highlighted by a more detailed account.
\end{abstract}



\section{Introduction: the background and around}  
\vspace*{-2mm}
The passage from commutative to noncommutative structures is a staple
in frontier physics and mathematics. The advent of quantum mechanics 
is of course an important example but it is now becoming increasingly 
clear that the appearance, 25 years ago \cite{FLS76, BFFLS}, of 
deformation quantization, followed by numerous developments and 
metamorphoses, has been a major factor in the present trend. 
The epistemological background can be traced back \cite{Fl82,DS} to the 
deformation philosophy which Mosh\'e Flato developed in the early seventies, 
motivated by deep (decades old) physical ideas and (then recent) mathematical 
developments. Since that time, inspiring many, he has consistently pursued 
and promoted physical and mathematical consequences of that idea in 
several directions, of which deformation quantization is presently the
most widely recognized.  

Sir Michael Atiyah said recently, after Oscar Wilde, that mathematics and
physics are two communities separated by a common language. Sir Michael 
added that the two communities managed to communicate rather well until 
the beginning of the twentieth century, then became so separated that, 
half a century ago, Eugene Wigner marveled at what he called the 
unreasonable effectiveness of mathematics in physics. In the past decade 
however, in spite of the ever increasing ``Babel tower effect" in Science, 
some form of communication was developed, with a converse phenomenon: what 
may seem an unreasonable effectiveness of physics in mathematics, including
in such abstract fields as algebraic geometry. Now, if we remember that 
mathematics arose as an abstraction of our understanding of the physical 
world, neither effectiveness should be unreasonable. We shall see here that
deformation quantization is a perfect example of both of them. 

Physical theories have their domain of applicability defined by the relevant
distances, velocities, energies, etc. concerned. But the passage from 
one domain (of distances, etc.) to another does not happen in an uncontrolled 
way. Rather, experimental phenomena appear that cause a paradox and 
contradict accepted theories. Eventually a new fundamental constant enters 
and the formalism is modified. Then the attached structures (symmetries, 
observables, states, etc.) \textsl{deform} the initial structure. Namely, 
we have a new structure which in the limit, when the new parameter goes to 
zero, coincides with the previous formalism. The only question is, in which 
category do we seek for deformations? Usually physics is rather 
conservative and if we start e.g. with the category of associative 
or Lie algebras, we tend to deform in the same category, but there are
important examples of generalizations of this principle (e.g. quantum groups 
are deformations of Hopf algebras).

The discovery of the non-flat nature of Earth may be the first example 
of this phenomenon. Closer to us, the paradox coming from the Michelson 
and Morley experiment (1887) was resolved in 1905 by Einstein with the 
special theory of relativity: there the Galilean geometrical 
symmetry group of Newtonian mechanics is deformed to the Poincar\'e group, 
the new fundamental constant being $c^{-1}$ where $c$ is the velocity of 
light in vacuum. At around the same time the Riemann surface theory can be 
considered as a first mathematical example of deformations, even if 
deformations became systematically studied in the mathematical literature 
only at the end of the fifties with the profound works of Kodaira and 
Spencer \cite{KS} on deformations of complex analytic structures. 
Now, when one has an action on a geometrical structure, it is natural to 
try and ``linearize" it by inducing from it an action on an algebra of 
functions on that structure. This is implicitly what Gerstenhaber did 
in 1963 \cite{Ge} with his definition and thorough study of deformations 
of rings and algebras. It is in the Gerstenhaber sense that the Galileo 
group is deformed to the Poincar\'e group; that operation is the inverse 
of the notion of group contraction introduced ten years before, empirically, 
by \'In\"on\"u and Wigner \cite{IW}. This fact triggered strong interest 
for deformation theory in France among a number of theoretical physicists, 
including Flato who had just arrived from the Racah school and knew well
the effectiveness of symmetry in physical problems.  

In 1900, as a last resort to explain the black body radiation, Planck 
proposed the quantum hypothesis: the energy of light is not emitted 
continuously but in quanta proportional to its frequency. He wrote $h$ for
the proportionality constant which bears his name. This paradoxical
situation got a beginning of a theoretical basis when, again in 1905, 
Einstein came with the theory of the photoelectric effect. Around 1920, 
an ``agr\'eg\'e d'histoire'', Prince Louis de Broglie, was introduced to 
the photoelectric effect, together with the Planck--Einstein relations 
and the theory of relativity, in the laboratory of his much older brother, 
Maurice duc de Broglie. This led him, in 1923, to his discovery of the 
duality of waves and particles, which he described in his celebrated 
Thesis published in 1925, and to what he called `m\'ecanique ondulatoire'. 
German and Austrian physicists, in particular, Hermann Weyl, Werner 
Heisenberg and Erwin Schr\"odinger, transformed it into the quantum 
mechanics that we know, where the observables are operators in Hilbert 
spaces of wave functions.

Intuitively, classical mechanics is the limit of quantum mechanics 
when $\hbar=\frac{h}{2\pi}$ goes to zero. But how can this be realized 
when in classical mechanics the observables are functions over phase space 
(a Poisson manifold) and not operators? The deformation philosophy promoted 
by Flato showed the way: one has to look for deformations of algebras of 
functions over Poisson manifolds endowed with the Poisson bracket and 
realize, in an \textsl{autonomous} manner, quantum mechanics there. This 
required, as a preliminary, a detailed study of the corresponding cohomology 
spaces. As a first step, in 1974, the cochains were assumed \cite{FLS} to 
be 1-differentiable (given by bidifferential operators of order $(1,1)$). 
This fell short of the solution but inspired Vey \cite{Ve} who was able in 
1975 (for symplectic manifolds with vanishing third Betti number) to show the 
existence of such differentiable deformations. Doing so, he rediscovered
a formula for the deformed bracket (the sine of the Poisson bracket) that 
he did not know had been obtained in an entirely different context 
by Moyal \cite{Mo} in 1949. The technical obstacle (a solution of which, 
in an algebraic context, could later be traced back to a result hidden in 
a fundamental paper \cite{HKR} from 1962) was lifted and deformation 
quantization could be developed \cite{FLS76,BFFLS} in 1976--78. In this 
approach, as we shall see in the next section, \textsl{quantization is 
a deformation} of the associative (and commutative) product of classical 
observables (functions on phase space) driven by the Poisson bracket, 
namely a star-product.

In a context which then seemed unrelated to deformation theory, 
pseudodifferential operators were introduced  also at the end of the fifties
and became a very hot subject in mathematics thanks to the publication in 
1963 of the first index theorems by Atiyah and Singer \cite{AS}, which 
express an analytically defined index in topological terms. The composition 
of symbols of pseudodifferential operators, an important ingredient 
in the proof, is a nontrivial example of a star-product, but this fact 
was noticed only years later, after deformation quantization was 
introduced. 

There have been many generalizations of the original index theorems, 
including algebraic versions developed in particular by Connes \cite{Co} 
in the context of noncommutative geometry, a natural continuation of 
his important works of the seventies on operator algebras that were 
motivated by physical problems. It was developed shortly after the 
appearance of deformation quantization, using cyclic rather than plain 
Hochschild cohomology. If one adds to star-products the notion of trace
(see Section 3.2.2), which in this case comes from integration over the 
manifold on which the functions are defined, one gets closed 
star-products \cite{CFS} (classified by cyclic cohomology) that 
provide \cite{Co} other examples of algebras falling in the framework 
of noncommutative geometry. Now the Gelfand isomorphism provides 
a realization of a commutative algebra as an algebra of functions 
over a manifold, i.e., its spectrum. A natural but difficult question is 
then to develop a theory of noncommutative manifolds. Truly nontrivial 
examples (if we exclude the group case, cf. below)
are appearing only recently (see e.g. \cite{Ma,CDV}).  

Simple Lie groups and algebras are rigid for the Gerstenhaber notion of 
deformation but if one goes to the category of Hopf algebras, they can be 
deformed. This is what Drinfeld \cite{Dr} realized with a notion to
which he gave the spectacular (albeit somewhat misleading) name of 
quantum group, by considering star-products deforming the product 
in a Hopf algebra of functions on Lie groups having a compatible Poisson 
structure. The dual approach of deforming the coproduct in a closure of 
enveloping algebras, generalizing the $\mathfrak{su}(2)$ example 
discovered empirically \cite{KR} in 1981 in relation with quantum 
inverse scattering problems, was also taken independently by Jimbo \cite{Ji}.
The domain has since known an extensive development, at the beginning 
mostly in Russia \cite{FRT}. Numerous applications to physics were developed 
in many fields and several metamorphoses occurred both in the mathematical 
methods used and in the concept itself. The subject now covers thousands of 
references, a fraction of which can be found in \cite{ES,Maj,SS}, which are
recent complementary expositions (\cite{SS} has 1264 references!) We shall not 
attempt to develop further this important avatar of deformation quantization.

\section{Genesis and first developments}   
\subsection{A few preliminaries}
In this part we recall some very classical notions, definitions and 
properties, for the sake of completeness and the benefit of theoretical 
physicists who might not be familiar with them but would be willing to 
take up the mathematical jargon required without having to read an extensive
literature. It may be skipped by many readers. As for all of this paper, 
the interested reader is strongly encouraged to study the references,
and the references quoted in references (etc.), at least in the topics that 
are closest to his interests: the paper is meant to be mostly a starting 
point for newcomers in the domains covered.
\vspace*{3mm}

\noindent\textbf{2.1.1 \ Hochschild and Chevalley--Eilenberg cohomologies.}
When $A$ is an \textsl{associative} algebra (over some commutative ring
${\mathbb K}$), we can consider it as a module over itself with
the adjoint action (algebra multiplication) and we shall here define 
cohomology in that context. The generalization to cohomology valued 
in a general module is straightforward.
A Hochschild \textsf{$p$-cochain} is a $p$-linear map $C$ from $A^p$ to  
$A$ and its \textsf{coboundary} $bC$ is a $(p+1)$-cochain given by
\begin{eqnarray} \label{bH}
bC(u_0,\ldots,u_p)&=&u_0C(u_1,\ldots,u_p)-C(u_0u_1,u_2,\ldots,u_p)+
\cdots\nonumber\\
&+&\hspace*{-2mm}(-1)^pC(u_0,u_1,\ldots,u_{p-1}u_p)+(-1)^{p+1}C(u_0,\ldots,
u_{p-1})u_p.
\end{eqnarray}
One checks that we have here what is called a complex, i.e. $b^2=0$.
We say that a $p$-cochain $C$ is a \textsf{$p$-cocycle} if $bC=0$. We denote
by $Z^p(A,A)$ the space of $p$-cocycles and by $B^p(A,A)$ the
space of those $p$-cocycles which are coboundaries (of a $(p-1)$-cochain).
The $p$th \textsf{Hochschild cohomology} space (of $A$ valued in $A$ seen as
a bimodule) is defined as $H^p(A,A)=Z^p(A,A)/B^p(A,A)$. \textsf{Cyclic
cohomology} is defined using a bicomplex which includes the Hochschild
complex. We refer to \cite{Co} for a detailed treatment.

For \textsl{Lie algebras} (with bracket $\{\cdot,\cdot\}$) one has a similar
definition, due to Chevalley and Eilenberg \cite{CE}. The (ad-valued) 
$p$-cochains are here skew-symmetric, i.e. linear maps 
$B: \wedge^p A\longrightarrow A$, and the Chevalley coboundary operator 
$\partial$ is defined on a $p$-cochain $B$ by (where ${\hat{u}_j}$ means 
that $u_j$ has to be omitted):
\begin{eqnarray}
\partial C(u_0,\ldots,u_p)=
&\sum_{j=0}^p&(-1)^j\{u_j,C(u_0,\ldots,{\hat{u}_j},\ldots,u_p)\}\nonumber\\
&+&\sum_{i<j}(-1)^{i+j}C(\{u_i,u_j\},u_0,\ldots,{\hat{u}_i},\ldots,
{\hat{u}_j},\ldots,u_p).
\end{eqnarray}
Again one has a complex ($\partial^2=0$), cocycles and coboundaries spaces
$\mathcal{Z}^p$ and $\mathcal{B}^p$ (resp.) and by quotient the 
\textsf{Chevalley cohomology} spaces $\mathcal{H}^p(A,A)$, or in short 
$\mathcal{H}^p(A)$; the collection of all cohomology spaces will be denoted 
here $\mathcal{H}^\bullet$, or $H^\bullet$ for the Hochschild cohomology.
\vspace*{3mm}

\noindent\textbf{2.1.2 \ Gerstenhaber theory of deformations of algebras.}
Let $A$ be an algebra. By this we mean an \textsl{associative, Lie} or 
\textsl{Hopf} algebra, or a \textsl{bialgebra} [an associative algebra $A$ 
where one has in addition a coproduct $\Delta : A \longrightarrow A \otimes A$ 
and the obvious compatibility relations, see e.g. \cite{ES} for a precise
definition]. Whenever needed we assume it is also a \textsl{topological} 
algebra, i.e., endowed with a locally convex topology for which all 
needed algebraic laws are continuous. For simplicity we may think
that the base (commutative) ring ${\mathbb K}$ is the field of complex
numbers ${\mathbb C}$ or that of the real numbers ${\mathbb R}$. Extending it
to the ring ${\mathbb K}[[\lambda]]$ of formal series in some parameter 
$\lambda$ gives the module ${\tilde A} = A[[\lambda]]$, on which we can 
consider the preceding various algebraic (and topological) structures. 
In a number of instances we also need to look at $A[\lambda^{-1},\lambda]]$, 
formal series in $\lambda$ and polynomials in $\lambda^{-1}$ (considered 
at first as an independent parameter). 
\vspace{2mm}

\noindent\textsl{2.1.2.1. \ Deformations and cohomologies}.
A concise formulation of a Gerstenhaber deformation of an algebra (which
we shall call in short a \textsf{DrG-deformation} whenever a confusion may
arise with more general deformations) is \cite{Ge,GS,BFGP}:
\begin{definition}\label{defidefo}
A \textsf{deformation} of such an algebra $A$ is a 
${\mathbb K}[[\lambda]]$-algebra
${\tilde A}$ such that ${\tilde A}/\lambda {\tilde A} \approx A$.
Two deformations ${\tilde A}$ and ${\tilde A'}$ are said \textsf{equivalent} 
if they are isomorphic over ${\mathbb K}[[\lambda]]$ and ${\tilde A}$ is said 
\textsf{trivial} if it is isomorphic to the original algebra $A$ considered 
by base field extension as a ${\mathbb K}[[\lambda]]$-algebra.
\end{definition}
Whenever we consider a topology on $A$, ${\tilde A}$ is supposed to be
topologically free. For associative (resp. Lie) algebras, Definition 
\ref{defidefo} tells us that there exists a new product $\star$ (resp. 
bracket $[\cdot,\cdot]$) such that the new (deformed) algebra is again 
associative (resp. Lie). Denoting the original composition laws by ordinary 
product (resp. $\{\cdot,\cdot\}$) this means that, for $u,v\in A$ we have:
\begin{eqnarray}
u\star v &=& uv + \sum_{r=1}^\infty \lambda^r C_r(u,v) \label{a}\\
\left[u, v\right] &=& \{u,v\} + \sum_{r=1}^\infty \lambda^r B_r(u,v)\label{l}
\end{eqnarray}
the $C_r$ being Hochschild 2-cochains and the $B_r$ (skew-symmetric)
Chevalley 2-cochains, such that for $u,v,w\in A$ we have
$(u\star v) \star w=u\star (v\star w)$ and ${\cal{S}}[[u,v],w]=0$, where
${\cal{S}}$ denotes summation over cyclic permutations (we extend (\ref{a}) 
and (\ref{l}) to $A[[\lambda]]$ by ${\mathbb K}[[\lambda]]$-linearity). 
At each level $r$ we therefore need to fulfill the equations ($j,k\geq1$):
\begin{eqnarray}
D_r(u,v,w)&\equiv&\sum_{j+k=r}(C_j(C_k(u,v),w)-C_j(u,C_k(v,w))=bC_r(u,v,w)
\label{adefcond}\\
E_r(u,v,w)&\equiv&\sum_{j+k=r}{\cal{S}}B_j(B_k(u,v),w)=\partial B_r(u,v,w)
\label{ldefcond}
\end{eqnarray}
where $b$ and $\partial$ denote (respectively) the Hochschild and Chevalley
coboundary operators. In particular we see that for $r=1$ the driver
$C_1$ (resp. $B_1$) must be a 2-cocycle. Furthermore, assuming one has shown
that (\ref{adefcond}) or (\ref{ldefcond}) are satisfied up to some order $r=t$,
a simple calculation shows that the left-hand sides for $r=t+1$ are then
3-cocycles, depending only on the cochains $C_k$ (resp. $B_k$) of order
$k\leq t$. If we want to extend the deformation up to order $r=t+1$
(i.e. to find the required 2-cochains $C_{t+1}$ or $B_{t+1}$), this cocycle
has to be a coboundary (the coboundary of the required cochain):
\textsl{The obstructions to extend a deformation from one step to the next lie
in the 3-cohomology}. In particular (and this was Vey's trick) if one can
manage to pass always through the null class in the 3-cohomology, a cocycle
can be the driver of a full-fledged (formal) deformation.

For a (topological) \textsl{bialgebra}, denoting by $\otimes_\lambda$ the 
tensor product of ${\mathbb K}[[\lambda]]$-modules, we can identify
${\tilde A}\, {\hat{\otimes}}_\lambda {\tilde A}$ with
$(A\, {\hat{\otimes}}A)[[\lambda]]$, where ${\hat{\otimes}}$ denotes the 
algebraic tensor product completed with respect to some operator topology 
(e.g. projective for Fr\'echet nuclear topology), we similarly have a deformed
coproduct ${\tilde \Delta } = \Delta + \sum_{r=1}^\infty \lambda^r D_r$,
$D_r \in {\cal L}(A, A {\hat{\otimes}}A)$ and in this context appropriate
cohomologies can be introduced. As we have said we shall not elaborate on 
these, nor on the additional requirements for Hopf algebras, referring for 
more details to original papers and books.
\vspace{2mm}

\noindent\textsl{2.1.2.2.\ Equivalence} means that there is an isomorphism
$T_\lambda=I+\sum_{r=1}^\infty \lambda^r T_r$, $T_r\in{\cal L}(A,A)$ so that
$T_\lambda(u\star'v)=(T_\lambda u\star T_\lambda v)$ in the associative case, 
denoting by $\star$ (resp. $\star'$) the deformed laws in ${\tilde A}$ 
(resp. ${\tilde A'}$); and similarly in the Lie case.
In particular we see (for $r=1$) that a deformation is trivial at order 1 if
it starts with a 2-cocycle which is a 2-coboundary. More generally, exactly
as above, we can show \cite{BFFLS} that if two deformations are equivalent
up to some order $t$, the condition to extend the equivalence one step
further is that a 2-cocycle (defined using the $T_k$, $k\leq t$) is the
coboundary of the required $T_{t+1}$ and therefore \textsl{the obstructions to
equivalence lie in the 2-cohomology}. In particular, if that space is null,
all deformations are trivial.

An important property is that a \textsl{deformation of an associative 
algebra with unit} (what is called a unital algebra) is again unital, 
and \textsl{equivalent to a deformation with the same unit}.
This follows from a more general result of Gerstenhaber (for deformations
leaving unchanged a subalgebra); a proof can be found in \cite{GS}.

In the case of (topological) \textsl{bialgebras} or \textsl{Hopf}
algebras, \textsl{equivalence} of deformations has to be understood as
an isomorphism of (topological) ${\mathbb K}[[\lambda]]$-algebras, the 
isomorphism starting with the identity for the degree 0 in $\lambda$. 
A deformation is again said \textsf{trivial} if it is equivalent to that 
obtained by base field extension. For Hopf algebras the deformed algebras 
may be taken (by equivalence) to have the same unit and counit, but in 
general not the same antipode.
\vspace{2mm}

\noindent\textsl{2.1.2.3.\ Homotopy of deformations}. Recently Gerstenhaber 
has considered (for reasons that are related to the so-called Donald--Flanigan 
conjecture, see \cite{GG,GGS}) the question of (formal) 
\textsf{compatibility} of deformations, a kind of homotopy in the variety 
of algebras between two deformations (\ref{a}) with parameters 
$\lambda$ and $\lambda'$ and cochains $C_r$ and $C'_r$. 
By this he means a 2--parameter deformation of the form 
\begin{equation} \label{2param}
u\tilde{\star}v=uv+\lambda C_1(u,v) +\lambda' C'_1(u,v) + 
\sum_{r=2}^\infty \Phi_r(u,v;\lambda,\lambda')
\end{equation}
where each $\Phi_r$ is a polynomial of total degree $r$ in $\lambda$ and 
$\lambda'$, which reduces to the first one-parameter deformation when 
$\lambda'=0$ and to the second when $\lambda=0$.
At the first order the condition for this to hold (e.g. for associative 
algebras) is that the Gerstenhaber bracket \cite{Ge} $[C_1,C'_1]_G$ is a 
3-coboundary, and here also there are higher obstructions. As an example, 
it follows from \cite{HKR} that the Weyl algebra and the quantum plane 
are formally (but non analytically \cite{GZ}) compatible nonequivalent 
deformations of the polynomial algebra ${\mathbb{C}}[x,y]$. In (2.2.3.4) 
below we shall see another appearance of such a 2-parameter deformation
in a physical context.
\vspace*{3mm}

\noindent\textbf{2.1.3 \ The differentiable case, Poisson manifolds.}
Consider the algebra $N=C^\infty(X)$ of functions on a differentiable
manifold $X$. When we look at it as an associative algebra acting on itself
by pointwise multiplication, we can define the corresponding Hochschild
cohomologies. Now let $\Lambda$ be a skew-symmetric contravariant two-tensor
(possibly degenerate) defined on $X$, satisfying $[\Lambda,\Lambda]_{SN}=0$
in the sense of the Schouten--Nijenhuis bracket 
(a definition of which, both intrinsic and in terms of local coordinates, 
can be found in \cite{BFFLS,FLS}; see also (\ref{schouten}) below).

Then the inner product of $\Lambda$ with the 2-form $du\wedge dv$, 
$P(u,v)=i(\Lambda)(du\wedge dv)$, $u,v\in N$,  defines a
\textsf{Poisson bracket} $P$: it is obviously skew-symmetric, satisfies 
the Jacobi identity because $[\Lambda,\Lambda]_{SN}=0$ and the Leibniz
rule $P(uv,w) = P(u,w)v + uP(v,w)$. It is a bidifferential 2-cocycle for the
(general or differentiable) Hochschild cohomology of $N$, skewsymmetric
of order $(1,1)$, therefore nontrivial \cite{BFFLS} and thus defines an
infinitesimal deformation of the pointwise product on $N$. $(X,P)$ 
is called a \textsf{Poisson manifold} \cite{BFFLS,Li77}.

When $\Lambda$ is everywhere nondegenerate ($X$ is then necessarily of
even dimension $2\ell$), its inverse $\omega$ is a closed everywhere
nondegenerate 2-form ($d\omega=0$ is then equivalent to 
$[\Lambda,\Lambda]_{SN}=0$) and we say that $(X,\omega)$ is 
\textsf{symplectic}; $\omega^\ell$ is a volume element on $X$. 
Then one can in a consistent manner work with differentiable 
cocycles \cite{BFFLS,FLS} and the differentiable Hochschild 
$p$-cohomology space $H^p(N)$ is that of all 
skew-symmetric contravariant $p$-tensor fields \cite{HKR,Ve}, therefore is 
infinite-dimensional. Thus, except when $X$ is of dimension 2 (because 
then necessarily $H^3(N)=0$), the obstructions belong to an 
infinite-dimensional space where they may be difficult to trace. 
On the other hand, when $2\ell=2$, any 2-cocycle can be the driver of a 
deformation of the associative algebra~$N$: ``anything goes" in this case; 
some examples for ${\mathbb R}^2$ can be found already in \cite{Ve}.

Now endow $N$ with a Poisson bracket: we get a Lie algebra and can look at
its Chevalley cohomology spaces. Note that $P$ is bidifferential of order
$(1,1)$ so it is important to check whether the Gerstenhaber theory is
\textsl{consistent} when restricted to \textsl{differentiable} cochains (in 
both cases, of arbitrary order and order at most 1). The answer is positive; 
in brief, if a coboundary is differentiable, it is the coboundary of a 
differentiable cochain. Later it was found \cite{Co,Na} that assuming only
continuity gives the same type of results.

Since $P$ is of order (1,1), it was natural to study first the
1-differentiable cohomologies. When the cochains are restricted to be of
order $(1,1)$ with no constant term (then they annihilate constant functions,
which we write ``n.c." for ``null on constants") it was found \cite{Li}
that the Chevalley cohomology $\mathcal{H}^\bullet_{\rm{1-diff,n.c.}}(N)$
of the Lie algebra $N$ (acting on itself with the adjoint representation)
is exactly the de Rham cohomology $H^\bullet_{dR}(X)$. Thus
${\mathrm{dim}}\mathcal{H}^p_{\mathrm{1-diff,n.c.}}(N)=b_p(X)$, the $p$th 
Betti number of the manifold $X$. Without the n.c. condition one gets a 
slightly more complicated formula \cite{Li}; in particular if $X$ is 
symplectic with an exact 2-form $\omega=d\alpha$, one has here
$\mathcal{H}^p_{\rm{1-diff}}(N)=H^p_{dR}(X)\oplus H^{p-1}_{dR}(X)$.
\vspace{.1mm}

Then the ``three musketeers" \cite{FLS} could, in 1974, study  
\textsl{1-differentiable deformations} of the Poisson bracket Lie
algebra $N$ and develop some applications at the level of classical 
mechanics. In particular it was noticed that the ``pure" order $(1,1)$ 
deformations correspond to a deformation of the 2-tensor $\Lambda$; 
allowing constant terms and taking the deformed bracket in Hamilton equations 
instead of the original Poisson bracket gave a kind of friction term.

Shortly afterwards, triggered by that work, a ``fourth musketeer" J. Vey
\cite{Ve} noticed that in fact 
${\mathrm{dim}}\mathcal{H}^2_{\mathrm{diff,n.c.}}(N)=
1+{\mathrm{dim}}\mathcal{H}^2_{\mathrm{1-diff,n.c.}}(N)$
and that $\mathcal{H}^3_{\mathrm{diff,n.c.}}(N)$ is also finite-dimensional. 
He could then study more general differentiable deformations, rediscovering 
in the ${\mathbb R}^{2\ell}$ case the Moyal bracket. The latter was at that 
time rather ``exotic" and few authors (except for a number of physicists, 
see e.g.\cite{AW}) paid any attention to it. In \textsl{Mathematical Reviews} 
this bracket, for which \cite{Mo} is nowadays often quoted, is not even
mentioned in the review! Vey's work permitted to tackle the 
problem with differentiable deformations, and deformation quantization 
was born \cite{FLS76}.
\subsection{The founding papers and some follow-up}
\noindent\textbf{2.2.1 \ Classical mechanics and quantization.}
Classical mechanics, in Lagrangean or Hamiltonian form, assumed at first
implicitly a ``flat" phase space ${\mathbb R}^{2\ell}$, or at least considered
only an open connected set thereof. Eventually more general configurations
were needed and so the mathematical notion of manifold, on which mechanics
imposed some structure, was used. This has lead in particular to using the
notions of symplectic and later of Poisson manifolds, which were introduced 
also for purely mathematical reasons. One of these reasons has to do with 
families of infinite-dimensional Lie algebras, which date back
to works by \'Elie Cartan at the beginning of last century and regained a
lot of popularity (including in physics) in the past 30 years.

In 1927, Weyl came out with his quantization rule \cite{We}. If we start with 
a classical observable $u(p,q)$, some function on (flat) phase space 
${\mathbb R}^{2\ell}$ (with $p, q \in {\mathbb R}^{\ell}$), one can associate
to it an operator (the corresponding quantum observable) $\Omega(u)$ in the
Hilbert space $L^2({\mathbb R}^{\ell})$ by the following general recipe:
\begin{equation} \label{weyl}
 u \mapsto \Omega_w(u) = \int_{{\mathbb R}^{2\ell}} \tilde{u}(\xi,\eta)
{\exp}(i(P.\xi + Q.\eta)/\hbar)w(\xi,\eta)\ d^\ell \xi d^\ell \eta
\end{equation}
where $\tilde{u}$ is the inverse Fourier transform of $u$, $P_\alpha$ and
$Q_\alpha$ are operators satisfying the canonical commutation relations
$[P_\alpha , Q_\beta] = i\hbar\delta_{\alpha\beta}$
($\alpha, \beta = 1,...,\ell$), $w$ is a weight function and the integral
is taken in the weak operator topology.  What is now called normal ordering
corresponds to choosing the weight
$w(\xi,\eta) = {\exp}(-\frac{1}{4}(\xi^2 \pm \eta^2))$,
standard ordering (the case of the usual pseudodifferential operators in
mathematics) to the weight $w(\xi,\eta) = {\exp}(-\frac{i}{2}\xi\eta)$ and 
the original Weyl (symmetric) ordering to $w = 1$.
An inverse formula was found shortly afterwards by Eugene Wigner \cite{Wi}
and maps an operator into what mathematicians call its symbol by a kind
of trace formula. For example $\Omega_1$ defines an isomorphism of Hilbert
spaces between $L^2({\mathbb R}^{2\ell})$ and Hilbert--Schmidt operators on
$L^2({\mathbb R}^{\ell})$ with inverse given by
\begin{equation} \label{EPW}
u=(2\pi\hbar)^{--\ell}\, {\rm{Tr}}[\Omega_1(u)\exp((\xi.P+\eta.Q)/i\hbar)]
\end{equation}
and if $\Omega_1(u)$ is of trace class one has
${\rm{Tr}}(\Omega_1(u))=(2\pi\hbar)^{-\ell}\int u \, \omega^\ell
\equiv \mathrm{Tr}_\mathrm{\scriptstyle M}(u)$, the ``Moyal trace", where
$\omega^\ell$ is the (symplectic) volume $dx$ on ${\mathbb R}^{2\ell}$.
Numerous developments followed in the direction of phase-space methods, many
of which can be found described in \cite{AW}. Of particular interest to us
here is the question of finding a physical interpretation to the classical 
function $u$, symbol of the quantum operator $\Omega_1(u)$; this was the 
problem posed (around 15 years after \cite{Wi}) by Blackett to his student 
Moyal. The (somewhat naive) idea to interpret it as a probability density
had of course to be rejected (because $u$ has no reason to be positive)
but, looking for a direct expression for the symbol of a quantum commutator,
Moyal found \cite{Mo} what is now called the Moyal bracket:
\begin{equation}
M(u,v) = \lambda^{-1} \sinh(\lambda P)(u,v) =
P(u,v) + \sum^\infty_{r=1}\lambda^{2r}P^{2r+1} (u,v)   \label{Moyal}
\end{equation}
where $2\lambda=i\hbar$, $P^r(u,v)=\Lambda^{i_1j_1}\ldots
\Lambda^{i_rj_r}(\partial_{i_1\ldots i_r}u)(\partial_{j_1\ldots j_r} v)$
is the $r^{\rm{th}}$ power ($r\geq1$) of the Poisson bracket
bidifferential operator $P$, $i_k, j_k = 1,\ldots,2\ell$, $k=1,\ldots,r$
and $\Lambda = {0\;-I\choose I\;\;\; 0}$.
To fix ideas we may assume here $u,v\in C^\infty({\mathbb R}^{2\ell})$ and the
sum taken as a formal series (the definition and convergence for various
families of functions $u$ and $v$ was also studied, including in
\cite{BFFLS}). A similar formula for the symbol of a product
$\Omega_1(u)\Omega_1(v)$ had been found a little earlier \cite{Gr} and
can now be written more clearly as a (Moyal) \textsf{star-product}:
\begin{equation} \label{star}
u \star_M v = \exp(\lambda P)(u,v) = uv + 
\sum^\infty_{r=1}\lambda^{r}P^{r}(u,v).
\end{equation}
Several integral formulas for the star-product have been introduced (we shall
come back later to this question). In great part after deformation quantization
was developed, the Wigner image of various families of operators (including 
all bounded operators on $L^2({\mathbb R}^{\ell})$) was studied and 
an adaptation to Weyl ordering of the mathematical notion of
pseudodifferential operators (ordered, like differential operators, 
``first $q$, then $p$") had been developed. We shall not give here
more details on those aspects, but for completeness mention a related
development. 
Starting from field theory, where normal (Wick) ordering is essential
(the role of $q$ and $p$ above is played by $q\pm ip$), Berezin
\cite{Be,BS} developed in the mid-seventies an extensive study of what
he called ``quantization", based on the correspondence principle and Wick
symbols. It is essentially based on K\"ahler manifolds and related to
pseudodifferential operators in the complex domain \cite{BG}.
However in his approach, as in the studies of various orderings \cite{AW}, 
the important concepts of \textsl{deformation} and \textsl{autonomous} 
formulation of quantum mechanics in general phase space are absent.

Quantization involving more general phase spaces was treated, in a 
somewhat systematic manner, only with Dirac constraints \cite{Di}: second
class Dirac constraints restrict phase space from some ${\mathbb R}^{2\ell}$
to a symplectic manifold $W$ imbedded in it (with induced symplectic form),
while first class constraints further restrict to a Poisson manifold
with symplectic foliation (see e.g. \cite{FLS76}). The question of 
quantization on such manifolds was certainly treated by many authors 
(including \cite{Di}) but did not go beyond giving some (often useful) 
recipes and hoping for the best.

A first systematic attempt started around 1970 with what was called soon
afterwards \textsl{geometric quantization} \cite{Ko}, a by-product of Lie 
group representations theory where it gave significant results.
It turns out that it is geometric all right, but its scope as far as
quantization is concerned has been rather limited since few classical
observables could be quantized, except in situations which amount
essentially to the Weyl case considered above.
In a nutshell one considers phase-spaces $W$ which are coadjoint orbits of
some Lie groups (the Weyl case corresponds to the Heisenberg group with the
canonical commutation relations and ${\mathfrak{h}}_\ell$ as Lie algebra); 
there one defines a ``prequantization" on the Hilbert space $L^2(W)$ and tries
to halve the number of degrees of freedom by using polarizations (often
complex ones, which is not an innocent operation as far as physics is
concerned) to get a Lagrangean submanifold ${\cal L}$ of dimension half of
that of $W$ and quantized observables as operators in $L^2({\cal L})$.
A recent exposition can be found in \cite{Wo}.
\vspace*{3mm}

\noindent\textbf{2.2.2 \ Star-products.} 
Star-products as a deformation of the usual product of functions on a
phase-space meant for an understanding of quantum mechanics were 
introduced in \cite{FLS76} and their relevance examplified by the
founding papers \cite{BFFLS}, which included significant applications.
Their general definition is as follows.

Let $X$ be a differentiable manifold (of finite, or possibly infinite,
dimension). We assume given on $X$ a \textsf{Poisson structure} 
(a Poisson bracket $P$).
\begin{definition}
A \textsf{star-product} on $X$ is a deformation of the associative algebra of 
functions $N=C^\infty(X)$ of the form $\star = \sum_{n=0}^\infty \lambda^n C_n$
where $C_0(u,v)=uv$, $C_1(u,v)-C_1(v,u)=2P(u,v)$, $u,v\in N$, and the
$C_n$ are bidifferential operators (locally of finite order).
We say a star-product is \textsf{strongly closed} if
$\int_X (u\star v-v\star u) dx =0$ where $dx$ is a volume element on $X$.
\end{definition}
The parameter $\lambda$ of the deformation is taken to be 
$\lambda= {\frac{i}{2}}\hbar$ in physical applications.
\vspace{2mm}

\noindent \textsl{2.2.2.1.\ a}. Using equivalence one may take $C_1=P$. 
That is the case of Moyal, but other orderings like standard or normal 
do not verify this condition (only the skew-symmetric part of $C_1$ is $P$).
Again by equivalence, in view of Gerstenhaber's result on the unit (cf. 
2.1.2.2), we may take cochains $C_r$ which are without constant term 
(what we called n.c. or null on constants). In fact, in the original paper 
\cite{BFFLS}, only this case was considered and the accent was put on 
``Vey products", for which the cochains $C_r$ have the same parity as $r$ and
have $P^r$ for principal symbol in any Darboux chart, with $X$ symplectic.
\vspace{1mm}  

\textsl{b}. It is also possible to consider star-products for which the
cochains $C_n$ are allowed to be slightly more general. Allowing them
to be \textsf{local} ($C_n(u,v)=0$ on any open set where $u$ or $v$ vanish)
gives nothing new. Note that this is not the same as requiring
the whole associative product to be local (in fact the latter
condition is very restrictive and, like true pseudodifferential operators,
a star-product is a nonlocal operation). In some cases (e.g. for star
representations of Lie groups) it may be practical to consider
pseudodifferential cochains. As far as the cohomologies are concerned, 
as long as one requires at least continuity for the cochains, 
the theory is the same as in the differentiable case \cite{Na}. 

Also, due to formulas like (\ref{EPW}) and the relation with Lie algebras
(see below), it is sometimes convenient to take
${\mathbb K}[\lambda^{-1},\lambda]]$ (Laurent series in $\lambda$, polynomial 
in $\lambda^{-1}$ and formal series in $\lambda$) for the ring on which 
the deformation is defined. Again, this will not change the theory.
\vspace{1mm}

\textsl{c}. By taking the corresponding commutator
$[u,v]_\lambda=(2\lambda)^{-1}(u\star v - v \star u)$, since the 
skew-symmetric part of $C_1$ is $P$, we get a deformation of the Poisson 
bracket Lie algebra $(N,P)$. This is a crucial point because
(at least in the symplectic case) we know the needed Chevalley cohomologies
and (in contradistinction with the Hochschild cohomologies) they are small
\cite{Ve,DLG}. The interplay between both structures gives existence
and classification; in addition it will explain why (in the symplectic case)
the classification of star-products is based on the 1-differentiable
cohomologies, hence ultimately on the de Rham cohomology of the manifold.
\vspace{2mm}

\noindent \textsl{2.2.2.2. \ Invariance and covariance.}
The Poisson bracket $P$ is (by definition) invariant under all 
symplectomorphisms, i.e. transformations of the manifold
$X$ generated by the flows $x_u=i(\Lambda)(du)$ defined by Hamiltonians 
$u\in N$. But already on ${\mathbb R}^{2\ell}$ one sees easily that its 
powers $P^r$, hence also the Moyal bracket (\ref{Moyal}), are invariant 
only under flows generated by Hamiltonians $u$ which are polynomials of 
maximal order 2, forming the ``affine" symplectic Lie algebra
${\mathfrak{sp}}({\mathbb R}^{2\ell})\cdot{\mathfrak{h}_\ell}$.
For other orderings the invariance is even smaller 
(only ${\mathfrak{h}_\ell}$ remains). For general Vey products the first 
terms of a star-product are \cite{BFFLS}
$C_2 = P^2_\Gamma + bH$ and $C_3 = S^3_\Gamma + T +3\partial H$.
Here $H$ is a differential operator of maximal order 2, $T$ a 2-tensor
corresponding to a closed 2-form, $\partial$ the Chevalley coboundary
operator. $P^2_\Gamma$ is given (in canonical coordinates) by an expression
similar to $P^2$ in which usual derivatives are replaced by covariant
derivatives with respect to a given symplectic connection $\Gamma$
(a torsionless connection with totally skew-symmetric components when all
indices are lowered using $\Lambda$). $S^3_\Gamma$ is a very special cochain
given by an expression similar to $P^3$ in which the derivatives are
replaced by the relevant components of the Lie derivative of $\Gamma$ in
the direction of the vector field associated to the function ($u$ or $v$).
Fedosov's algorithmic construction \cite{Fe} shows that the symplectic
connection $\Gamma$ plays a role at all orders. Therefore the invariance
group of a star-product is a subgroup of the \textsl{finite-dimensional} group
of symplectomorphisms preserving a connection. Its Lie algebra is
${\mathfrak{g}}_{0}=\{a\in N; [a,u]_\lambda=P(a,u)\, \forall u\in N \}$, the
elements of which are \textsf{preferred observables}, i.e., Hamiltonians for 
which the classical and quantum evolutions coincide.
We are thus lead to look for a weaker notion and shall call a star-product
\textsf{covariant} under a Lie algebra ${\mathfrak{g}}$ of functions if
$[a,b]_\lambda = P (a,b)$ $\forall a,b \in {\mathfrak{g}}$. It can be shown
\cite{ACMP} that $\star$ is ${\mathfrak{g}}$-covariant iff there exists a
representation $\tau$ of the Lie group $G$ whose Lie algebra is
${\mathfrak{g}}$ into ${\rm{Aut}}(N[[\lambda]]; \star)$ such that
$\tau_g u = (Id_N + \sum^\infty_{r=1}\lambda^r\tau^r_g)(g.u)$ where
$g \in G, u \in N$, $G$ acts on $N$ by the natural action induced by the
vector fields associated with ${\mathfrak{g}}$, $(g\cdot u)(x)=u(g^{-1} x)$,
and where the $\tau^r_g$ are differential operators
on $W$. Invariance of course means that the geometric action preserves the
star-product: $g\cdot u \star g \cdot v = g\cdot (u \star v)$. This is the
basis for the theory of star representations which we shall briefly
present below, and of the relevance of star-products for important problems
in group theory.
\vspace*{3mm}

\noindent\textbf{2.2.3 \ Quantum mechanics.}
Let us start with a phase space $X$, a symplectic (or Poisson) manifold and
$N$ an algebra of classical observables (functions, possibly including
distributions if proper care is taken for the product). We shall call
\textsf{star quantization} a star-product on $X$ invariant (or sometimes only
covariant) under some Lie algebra ${\mathfrak{g}}_0$ of ``preferred 
observables". Invariance of the star-product ensures that the classical and
quantum evolutions of observables under a Hamiltonian 
$H \in {\mathfrak{g}}_0$ will coincide \cite{BFFLS}. The typical example
is the Moyal product on $W = {\mathbb{R}}^{2\ell}$.
\vspace{2mm}

\noindent \textsl{2.2.3.1 \ Spectrality.} Physicists want to get numbers
matching experimental results, e.g. for energy levels of a system. That is
usually achieved by describing the spectrum of a given Hamiltonian $\hat{H}$
supposed to be a self-adjoint operator so as to get a real spectrum and
so that the evolution operator (the exponential of $it \hat{H}$) is unitary 
(thus preserves probability amplitude). A similar spectral theory can
be done here, in an \textsl{autonomous manner}. The most efficient way to
achieve it is to consider \cite{BFFLS} the \textsf{star exponential}
(corresponding to the evolution operator)
\begin{equation}\label{starexp}
\mbox{Exp}(Ht) \equiv \sum^\infty_{n=0}
\frac{1}{n!}\left(\frac{t}{i\hbar}\right)^n(H\star)^n
\end{equation}
where $(H\star)^n$ means the $n^{\mbox{\footnotesize th}}$ star power of the
Hamiltonian $H \in N$ (or $N[[\lambda]]$).
Then one writes its Fourier--Stieltjes transform $d\mu$ (in the distribution
sense) as $\mbox{Exp}(Ht) = \int e^{\lambda t/i\hbar} d\mu (\lambda)$ and
defines \textsl{the spectrum of $(H /\hbar)$} as the \textsl{support $S$ of 
$d\mu$} (incidentally this is the definition given by L.~Schwartz for the 
spectrum of a distribution, out of motivations coming from Fourier analysis).
In the particular case when $H$ has discrete spectrum, the integral can
be written as a sum (see the top equation in (\ref{spec}) below for a
typical example): the distribution $d\mu$ is a sum of ``delta functions"
supported at the points of $S$ multiplied by the symbols of the
corresponding eigenprojectors.

In different orderings with various weight functions $w$ in (\ref{weyl}) one
gets in general different operators for the same classical observable $H$,
thus different spectra. For $X = {\mathbb{R}}^{2\ell}$ all orderings are
mathematically equivalent (to Moyal under the Fourier transform $T_w$
of the weight function $w$). This means that every observable $H$ will have
the same spectrum under Moyal ordering as $T_w H$ under the equivalent
ordering. But this does not imply physical equivalence, i.e. the fact that
$H$ will have the same spectrum under both orderings. In fact the opposite
is true: if two equivalent star-products are isospectral
(give the same spectrum for a large family of observables and all $\hbar$),
they are identical \cite{CFGS}.

It is worth mentioning that our definition of spectrum permits to define
a spectrum even for symbols of non-spectrable operators, such as the
derivative on a half-line which has different deficiency indices; this
corresponds to an infinite potential barrier (see also \cite{UUp}
for detailed studies of similar questions). That is one of the many
advantages of our autonomous approach to quantization.
\vspace{2mm}

\noindent \textsl{2.2.3.2.\ Applications.}  In quantum mechanics it is 
preferable to work (for $X = {\mathbb{R}}^{2\ell}$) with the star-product 
that has maximal symmetry, i.e. 
${\mathfrak{sp}}({\mathbb{R}}^{2\ell})\cdot{\mathfrak{h}}_\ell$
as algebra of preferred observables: the Moyal product. One indeed
finds \cite{BFFLS} that the star exponential of these observables
(polynomials of order $\leq 2$) is proportional to the usual exponential.
More precisely, if $H=\alpha p^2+\beta pq + \gamma q^2 \in {\mathfrak{sl}}(2)$
with $p,q \in  {\mathbb{R}}^{\ell}$, $\alpha,\beta, \gamma \in {\mathbb{R}}$,
setting $d = \alpha\gamma -\beta^2$ and $\delta = \vert d \vert^{1/2}$
one gets, summing the star exponential series and taking its Fourier 
decompositions (the sums and integrals appearing in the various expressions
being convergent as distributions, both in phase-space
variables and in $t$ or $\lambda$)
\begin{equation}\label{expalt}
\mbox{Exp}(Ht) = \left\{ \begin{array}{r@{\quad \mbox{for} \quad}l}
 (\cos\delta t)^{-l} \exp((H/i\hbar\delta)\tan(\delta t))& d > 0 
\vspace*{1mm}\\ \exp(Ht/i\hbar)& d = 0  \vspace*{1mm}\\
(\cosh\delta t)^{-l}\exp((H/i\hbar\delta)\tanh(\delta t))& d < 0
\end{array}\right.
\end{equation}
\begin{equation}\label{spec}
\mbox{Exp}(Ht) = \left\{\begin{array}{r@{\quad \mbox{for} \quad}l}
\sum^\infty_{n=0} \Pi^{(\ell)}_n e^{(n + \frac{\ell}{2})t} & d > 0 
\vspace*{2mm}\\
\int^\infty_{- \infty} e^{\lambda t/i\hbar} \Pi(\lambda, H)d\lambda & d < 0
\end{array}\right.
\end{equation}
We get the discrete spectrum $(n + \frac{\ell}{2})\hbar$ of the
\textsl{harmonic oscillator} and the continuous spectrum ${\mathbb{R}}$
for the dilation generator $pq$. The eigenprojectors $\Pi^{(\ell)}_n$ and
$\Pi (\lambda, H)$ are given \cite{BFFLS} by known special functions on
phase-space (generalized Laguerre and hypergeometric, multiplied by some
exponential). Formulas (\ref{expalt}) and (\ref{spec}) can, by analytic
continuation, be extended outside singularities and (as distributions) 
for singular values of $t$.

Other examples can be brought to this case by functional manipulations
\cite{BFFLS}. For instance the Casimir element $C$ of ${\mathfrak{so}}(\ell)$
representing \textsl{angular momentum}, which can be written
$C = p^2 q^2 - (pq)^2 - \ell(\ell - 1)\frac{\hbar^2}{4}$,
has $n (n + (\ell-2)) \hbar^2$ for spectrum. For the \textsl{hydrogen atom},
with Hamiltonian $H = \frac{1}{2} p^2 - \vert q \vert^{-1}$, the Moyal
product on ${\mathbb{R}}^{2\ell+2}$ ($\ell=3$ in the
physical case) induces a star-product on $X = T^*S^\ell$; the energy
levels, solutions of $(H-E) \star \phi = 0,$ are found from (\ref{spec})
and the preceding calculations for angular momentum to be (as they should,
with $\ell=3$) $E = \frac{1}{2} (n+1)^{-2} \hbar^{-2}$ for the discrete
spectrum, and $E \in {\mathbb{R}}^+$ for the continuous spectrum.

We thus have recovered, in a completely autonomous manner entirely within
deformation quantization, the results of ``conventional" quantum mechanics
in those typical examples (many more can be treated similarly).
It is worth noting that the term $\frac{\ell}{2}$ in the harmonic oscillator
spectrum, obvious source of divergences in the infinite-dimensional case,
disappears if the normal star-product is used instead of Moyal -- which
is one of the reasons it is preferred in field theory.
\vspace{2mm}

\noindent \textsl{2.2.3.3.\ Remark on convergence and the physical meaning 
of deformation quantization}. We have always considered
star-products as formal series and looked for convergence only in specific
examples, generally in the sense of distributions. The same applies
to star exponentials, as long as each coefficient in the formal series is
well defined. In the case of the harmonic oscillator or more generally
the preferred observables $H$ in Weyl ordering, the study is
facilitated by the fact that the powers $(H\star)^n$ are polynomials in $H$.

Integral formulas for star-products and star exponentials can 
be used to give a meaning to these expressions beyond the domains where
the sums are convergent or summable. \textsl{The sums must not be 
considered as perturbation series}, rather taken as a whole, whether
expressed as integral formulas or formal series. 

We stress that deformation quantization should not be seen as a mere 
reformulation of quantum mechanics or quantum theories in general. 
At the conceptual level, it is the true mathematical formulation of 
physical reality whenever quantum effects have to be taken into account. 
The above examples show that one can indeed perform important 
quantum mechanical calculations, in an autonomous manner, entirely within
deformation quantization -- and get the results obtained in conventional 
quantum mechanics. Whether one uses an operatorial formulation or some 
form of deformation quantization formulation is thus basically a 
practical question, which formulation is the most effective, at least 
in the cases where a Weyl or Wigner map exists. When such a map does not 
exist, a satisfactory operatorial formulation will be very hard to find, 
except locally on phase space, and deformation quantization is the solution. 

One can of course (and should in practical examples, as we have done here,
and also for algebraic varieties) look for small domains (in $N$) where one 
has convergence. We can then speak of ``strict" deformation quantization. 
In particular we can look for domains where pointwise
convergence can be proved; this was done e.g. for Hermitian symmetric spaces
\cite{CGR}. But it should be clearly understood that one can consider
wider classes of observables -- in fact, the latter tend to be physically more
interesting -- than those that fit in a strict $C^*$-algebraic approach. 
\vspace{2mm}

\noindent\textsl{2.2.3.4.\ A 2-parameter star-product and statistical 
mechanics.}
A relation similar to (\ref{2param}) had been considered for 
star-products in \cite{BFLS}, in connection with statistical mechanics.
Indeed, in view of our philosophy on deformations, a natural question to 
ask is their \textsf{stability}: Can deformations be further deformed, 
or does ``the buck stops there"? As we indicated at the beginning,
the answer to that question may depend on the context. Quantum groups are
an example, when dealing with Hopf algebras; properties of Harrison 
cohomology require going to noncontinuous cochains if one wants nontrivial
abelian deformations (see (4.3.3.3) and \cite{CF}). Here is a simpler example.

If one looks for deformations of the Poisson bracket Lie algebra $(N,P)$
one finds (under a mild assumption) that a further deformation of the Moyal 
bracket, with another deformation parameter $\rho$, is again a Moyal 
bracket for a $\rho$-deformed Poisson structure; in particular, 
for $X={\mathbb R}^{2\ell}$,
\textsl{quantum mechanics viewed as a deformation is unique and stable}.

Now, for the associative algebra $N$, the only \textsl{local} associative
composition laws are of the form $(u,v)\mapsto ufv$ for some $f\in N$. 
If we take $f=f_\beta\in N[[\beta]]$ we get a 0-differentiable
deformation (with parameter $\beta$) of the usual product. We were thus lead
\cite{BFLS} to look, starting from a product $\star_\lambda$, for a new
composition law $(u,v)\mapsto u{\tilde{\star}}_{\lambda,\beta}v=
u\star_\lambda f_{\lambda,\beta}\star_\lambda v$ with $f_{\lambda,\beta}=
\sum_{r=0}^\infty \lambda^{2r}f_{2r,\beta}\in N[[N[[\lambda^2]],\beta]]$,
where $f_{0,\beta}\equiv f_\beta \neq 0$ and $f_0=1$. The transformation
$u\mapsto T_{\lambda,\beta}u=f_{\lambda,\beta}\star_\lambda u$ intertwines 
$\star_\lambda$ and ${\tilde{\star}}_{\lambda,\beta}$ but it is not 
an equivalence of star-products because ${\tilde{\star}}_{\lambda,\beta}$
is not a star-product: it is a ($\lambda,\beta$)-deformation of the usual 
product with at first order in $\lambda$ the driver given by 
$P_\beta(u,v)=f_\beta P(u,v)+uP(f_\beta,v)-P(f_\beta,u)v$, 
a conformal Poisson bracket associated with a \textsf{conformal symplectic
structure} given by the 2-tensor $\Lambda_\beta=f_\beta\Lambda$ and the
vector $E_\beta=[\Lambda,f_\beta]_{SN}$.

We then start with a star-product, denoted $\star$, on some algebra ${\cal{A}}$ 
of observables, and take $f_{\lambda,\beta}=\exp_\star(c\beta H)$ with 
$c=-{\frac{1}{2}}$. The star exponential $\mbox{Exp}(Ht)$ defines an 
automorphism $u\mapsto\alpha_t(u)=\mbox{Exp}(-Ht)\star u\star\mbox{Exp}(Ht)$.
A \textsf{KMS state} $\sigma$ on ${\cal{A}}$ is a state (linear functional)
satisfying, $\forall a,b\in {\cal{A}}$, the Kubo--Martin--Schwinger condition
$\sigma(\alpha_t(a)\star b)=\sigma(b\star\alpha_{t+i\hbar\beta}(a))$.
Then the (quantum) KMS condition can be written \cite{BFLS}, 
with $[a,b]_\beta= (i/\hbar)(a{\tilde{\star}_{\lambda,\beta}}b
-b{\tilde{\star}}_{\lambda,\beta}a)$, 
simply $\sigma(g_{-\beta}\star[a,b]_\beta)=0$: up to a conformal factor, 
a KMS state is like a trace with respect to this new product. 
The (static) classical KMS condition is the limit for $\hbar=0$ of the 
quantum one. So we can recover known features of statistical mechanics 
by introducing a new deformation parameter $\beta=(kT)^{-1}$ and the 
related conformal symplectic structure. This procedure commutes with usual 
deformation quantization. The question was recently treated from
a more conventional point of view in deformation quantization (see e.g. 
\cite{BRW,Ws97}) using the notion of formal trace \cite{NT95}.
\vspace*{3mm}

\noindent\textbf{2.2.4 \ Star representations, a short overview.}
Let $G$ be a Lie group (connected and simply connected), acting by
symplectomorphisms on a symplectic manifold $X$ (e.g. coadjoint orbits
in the dual of the Lie algebra ${\mathfrak{g}}$ of $G$). The elements
$x,y \in {\mathfrak{g}}$ will be supposed to be realized by functions
$u_x, u_y$ in $N$ so that their Lie bracket $[x,y]_{\mathfrak g}$ is 
realized by $P(u_x, u_y)$. Now take a $G$-covariant star-product $\star$, 
that is $P(u_x,u_y) = [u_x, u_y]\equiv (u\star v-v\star u)/2\lambda$, 
which shows that the map ${\mathfrak{g}}\ni x \mapsto (2\lambda)^{-1}u_x\in N$
is a Lie algebra morphism. The appearance of $\lambda^{-1}$ here and in 
the trace (see 2.2.1) cannot be avoided and explains why we have often 
to take into account both $\lambda$ and $\lambda^{-1}$. 
We can now define the {\it star exponential}
\begin{equation}\label{stargp}
E (e^x) = {\mbox{Exp}}(x) \equiv
\sum^\infty_{n=0} (n!)^{-1} (u_x /2 \lambda)^{\star n}
\end{equation}
where $x \in {\mathfrak{g}}$, $e^x \in G$ and the power $\star n$ denotes the
$n^{\mbox{\footnotesize th}}$ star-power of the corresponding function. By the
Campbell--Hausdorff formula one can extend $E$ to a {\it group homomorphism}
$E : G \rightarrow (N[[\lambda,\lambda^{-1}]], \star)$ where, in the formal 
series, $\lambda$ and $\lambda^{-1}$ are treated as independent parameters 
for the time being. Alternatively, the values of $E$ can be taken in the 
algebra $({\cal P} [[\lambda^{-1}]], \star)$, where ${\cal P}$ is the 
algebra generated by ${\mathfrak g}$ with the $\star$-product 
(it is a representation of the enveloping algebra).

\noindent A \textsf{star representation} \cite{BFFLS} of $G$ is a distribution
${\cal E}$ (valued in ${\mathrm{Im}} E$) on $X$ defined by
\begin{displaymath}
D \ni f \mapsto {\cal E} (f) = \int_G f(g) E(g^{-1}) dg
\end{displaymath}
where $D$ is some space of test-functions on $G$. The corresponding
\textsf{character} $\chi$ is the (scalar-valued) distribution defined by
$D \ni f \mapsto \chi (f) = \int_X {\cal E} (f)d\mu $, $d\mu$ being a
quasi-invariant measure on $X$.

The character is one of the tools which permit a comparison with usual
representation theory. For semi-simple groups it is singular at the origin
in irreducible representations, which may require caution in computing the
star exponential (\ref{stargp}). In the case of the harmonic oscillator that
difficulty was masked by the fact that the corresponding representation of
${\mathfrak{sl}}(2)$ generated by $(p^2,q^2,pq)$ is integrable to a double
covering of ${\rm{SL}}(2,{\mathbb{R}})$ and decomposes into a sum
$D({\frac{1}{4}})\oplus D({\frac{3}{4}})$: the singularities at the origin
cancel each other for the two components.

This theory is now very developed, and parallels in many ways the usual
(operatorial) representation theory. A detailed account of all the results
would take us too far, but among the most notable one may quote:
\begin{itemize} \item[i)] 
An exhaustive treatment of \textsl{nilpotent} or \textsl{solvable exponential}
\cite{AC} and even \textsl{general solvable} Lie groups \cite{ACL}.
The coadjoint orbits are there symplectomorphic to ${\mathbb R}^{2\ell}$
and one can lift the Moyal product to the orbits in a way that is adapted
to the Plancherel formula. Polarizations are not required, and
``star-polarizations" can always be introduced to compare with usual theory.
Wavelets, important in signal analysis, are manifestations of
star-products on the (2-dimensional solvable) affine group of ${\mathbb{R}}$
or on a similar 3-dimensional solvable group \cite{BB}.
\item[ii)] 
%
For \textsl{semi-simple} Lie groups an array of results is already
available, including \cite{ACG,Mr} a complete treatment of the
\textsl{holomorphic discrete series} (this includes the case of compact
Lie groups) using a kind of Berezin dequantization, and scattered results
for specific examples. Similar techniques have also been used 
\cite{CGR,Ka} to find invariant star-products on K\"ahler and 
Hermitian symmetric spaces (convergent for an appropriate dense subalgebra). 
Note however, as shown by recent developments of unitary representations
theory (see e.g. \cite{Sch}),
that for semi-simple groups the coadjoint orbits alone are no more sufficient
for the unitary dual and one needs far more elaborate constructions.
\item[iii)]
%
For semi-direct products, and in particular for the Poincar\'e and
Euclidean groups, an autonomous theory has also been developed (see e.g.
\cite{ACM}).
\end{itemize}

Comparison with the usual results of ``operatorial" theory of Lie group
representations can be performed in several ways, in particular by
constructing an invariant Weyl transform generalizing (\ref{weyl}),
finding ``star-polarizations" that always exist, in contradistinction
with the geometric quantization approach (where at best one can
find complex polarizations), study of spectra (of elements in the center
of the enveloping algebra and of compact generators) in the sense of
(2.2.3.1), comparison of characters, etc. Note also in this context that the
pseudodifferential analysis and (non autonomous) connection with
quantization developed extensively by Unterberger, first in
the case of ${\mathbb{R}}^{2\ell}$, has been extended to the above
invariant context \cite{UU,UUp}. But our main insistence is that
the theory of star representations is an \textsl{autonomous} one that can be
formulated completely within this framework, based on coadjoint orbits
(and some additional ingredients when required).
\section{Existence and classification: geometric approach} \label{ga} 
\addtocounter{equation}{-15}
The first proof of the existence of star-products on any (real) symplectic 
manifold was given by De~Wilde and Lecomte \cite{DL83} and
is based on the idea of gluing local Moyal star-products defined on Darboux
charts. Their proof is algebraic and cohomological in essence. 
A more geometrical proof was provided by the construction of 
Omori, Maeda, and Yoshioka \cite{OMY91}. 
They have introduced the notion of Weyl bundle as a 
locally trivial fiber bundle whose fibers are Weyl algebras. 
A Weyl algebra is a unital $\mathbb{C}$-algebra generated by 
$1,\lambda,X_1,\ldots,X_n,Y_1,\ldots,Y_n$ satisfying the commutation relations
$[X_i,Y_j]=-\delta_{ij}\lambda$ and all other commutators are trivial.
The main step was to show the existence of a Weyl bundle $W$ over any 
symplectic manifold $M$. This was done by defining a special class of sections 
(the Weyl functions), first locally and then globally by showing that the class
of Weyl functions is stable under different choices of local trivialization. 
Then, the composition of Weyl sections induces a star-product on $M$. 
The construction of a Weyl bundle shares some features with the 
De~Wilde--Lecomte approach: in both cases, derivations play an important 
role and the basic idea is still to glue together local objects defined 
on Darboux charts by using cohomological arguments.

In 1985, motivated by index theory and independently of the previous proofs, 
Fedosov \cite{Fe1} announced a purely geometrical construction, also based 
on Weyl algebras, of star-products on a symplectic manifold. 
Unfortunately, the announcement appeared in a local volume with a summary
in Doklady and has stayed almost unknown for several years until a 
complete version was published in an international journal \cite{Fe}. 
The beautiful construction of Fedosov does provide a new proof of existence, 
but above all it gives a much better understanding of deformation 
quantization that allowed many major developments that will be 
discussed below. Let us first briefly review the Fedosov construction
(see also e.g. \cite{Ws}). For a comprehensive treatment and much more, 
we recommend Fedosov's book \cite{Fe3}.

\subsection{The construction of Fedosov}\label{fedosov}

Let $(M,\omega)$ be a $2m$-dimensional symplectic manifold. 
Let $\nabla$ be a symplectic connection, i.e., a linear connection
without torsion and preserving the symplectic form $\omega$.
At each point $x\in M$, the symplectic form endows the tangent
space $T_xM$ with a (constant) symplectic structure and we can define
the Moyal star-product $\ast_\lambda$ associated to $\omega_x$ on $T_xM$.
The space of formal series in $\lambda$ with polynomial coefficients on $T_xM$,
endowed with $\ast_\lambda$, gives us the  Weyl algebra $W_x$ 
associated to $T_xM$. The algebras $W_x$, $x\in M$, can be smoothly patched
and we get a fiber bundle $W=\cup_{x\in M}W_x$ on $M$, called the Weyl 
bundle over $M$. The fiberwise Moyal star-product endows naturally the 
space of sections $\Gamma(W)$ with a structure of unital associative 
algebra. The center of $\Gamma(W)$ can be identified as a vector space with
$C^\infty(M)[[\lambda]]$. The basic idea of Fedosov's construction is to
build from $\nabla$ a flat connection~$D$ on $W$ such that the algebra of
horizontal sections for $D$ induces a star-product on $M$.
This ``infinitesimal'' approach bypasses the need of gluing together
star-products defined on (large) Darboux charts.

In a local Darboux chart at $x\in M$, denote by  $(y^1,\ldots,y^{2m})$
the coordinates on $T_x M$. A section of $W$ can be locally written as 
\begin{equation}
 a(x,y)=\sum_{k\geq0,|\alpha|\geq0} \lambda^k a_{k,\alpha}(x) y^\alpha,
\end{equation}
where $\alpha = (\alpha_1,\ldots,\alpha_{2m})$ is a multi-index.
We shall need the bundle 
$$ \Gamma(W)\otimes \Lambda = 
\bigoplus_{0\leq q\leq 2m} \Gamma(W)\otimes \Lambda^q $$ 
of differential forms on $M$ taking their values in $W$. 
A section of $\Gamma(W)\otimes \Lambda$ can be expressed locally as
\begin{equation}\label{axydx}
a(x,y,dx)=
\sum_{k\geq0}\sum_{|\alpha|,|\beta|\geq0} 
\lambda^k a_{k,\alpha,\beta}(x) y^\alpha (dx)^\beta,
\end{equation}
where $(dx)^\beta\equiv dx^{\beta_1}\wedge\cdots\wedge dx^{\beta_q}$.
The natural extension of the product in $\Gamma(W)$ obtained by taking
the Moyal star-product $\ast_\lambda$ between the $y^i$'s and the wedge product
between the $dx^i$'s, defines a product $\circ$ which makes 
$\Gamma(W)\otimes \Lambda$ an algebra. There is also a structure of graded
Lie algebra on $\Gamma(W)\otimes \Lambda$ given by the bracket
$[a,b]=a\circ b - (-1)^{pq}b\circ a$ for
$a\in\Gamma(W)\otimes \Lambda^p$ and $b\in\Gamma(W)\otimes \Lambda^q$. 
It is just the natural extension of the bracket in $\Gamma(W)$ for the Moyal
star-product $\ast_\lambda$ to the differential forms taking values in $W$.

The algebra $\Gamma(W)\otimes \Lambda$ is  filtered in the following way.
Define the degrees of $\lambda$ and $y^i$ by $\deg(\lambda)=2$ and  
$\deg(y^i)=1$. Then $\Gamma(W)\otimes \Lambda$ has a filtration with 
respect to the total degree $2k + |\alpha|$ of the terms 
appearing in (\ref{axydx}):
$
\Gamma(W)\otimes \Lambda \supset \Gamma_1(W)
\otimes \Lambda\supset \Gamma_2(W)\otimes \Lambda \supset\cdots
$
Notice that the center of $\Gamma(W)\otimes \Lambda$ for the bracket $[\ ,\ ]$
is the set of sections independent of the $y^i$'s. We shall denote by $\sigma$
the projection of $\Gamma(W)\otimes \Lambda$ onto its center.

The exterior differential $d$ on $M$ extends to an antiderivation $\delta$
of the algebra $\Gamma(W)\otimes \Lambda$ by the formula 
$\delta(a)=\sum_j dx^j\wedge \partial_j a$, where $\partial_j=
\partial/\partial_{y^j}$. The antiderivation $\delta$ reduces by one unit
the degree of the filtration. Fedosov has also introduced another operation 
$\delta^*$ on $\Gamma(W)\otimes \Lambda$. It is defined by the formula
$\delta^*(a)=\sum_j y^j {\bf i}({\partial_j})(a)$, where ${\bf i}$ denotes
the interior product. One can check that $\delta^2= \delta^{*\ 2}=0$ 
and $\delta\delta^* +\delta^*\delta= (p+q)\mathrm{Id}$ on homogeneous elements
of degrees $p$ in $y^i$ and $q$ in $dx^i$. There is a Hodge-like decomposition
for the sections (\ref{axydx}). Denote by $a_{00}$ the term of the section $a$
in (\ref{axydx}) which is independent of the $y^i$'s and the $dx^i$'s.
By defining $\delta^{-1}$
on $(p,q)$-homogeneous elements by $\delta^{-1}=0$ if $p+q=0$ and 
$\delta^{-1}=\delta^*/(p+q)$ if $p+q>0$, every section $a$ admits the 
decomposition $a= \delta\delta^{-1}(a)+ \delta^{-1}\delta(a)+ a_{00}$.

The next step is to introduce another connection on $W$. First, the symplectic
connection $\nabla$ defines  on $W$ a connection $\hat{\nabla}$:
$\Gamma(W)\otimes\Lambda^q\rightarrow\Gamma(W)\otimes \Lambda^{q+1}$,
by  $\hat{\nabla}_i = \sum_i dx^i\wedge \nabla_i$. Fedosov introduces 
the new connection on $W$ as follows: for an element 
$\gamma\in \Gamma(W)\otimes \Lambda^1$ such that
$\sigma(\gamma)=0$, define the connection 
$D=\hat{\nabla} - [\gamma, \cdot]/i\lambda$.
The curvature of $D$ is  $\Omega=R+\hat{\nabla}\gamma -\gamma^2/i\lambda$, 
where $R$ is the curvature of  $\hat{\nabla}$. $\Omega$ is called the 
\textsf{Weyl curvature} of $D$. The usual relation 
$D^2 = - [\Omega,\cdot]/i\lambda$ holds.

Here comes the essential ingredient of Fedosov's construction. The idea is
to make quantum corrections to $D$ in order to make it flat, i.e., $D^2=0$.
In other words, it consists to make  the Weyl curvature $\Omega$ central 
(Abelian connection). Fedosov showed the existence of a flat Abelian 
connection on $W$, normalized by the condition $\sigma(r)=0$, of 
the form $D= -\delta + \hat{\nabla} - [r,\cdot]/i\lambda$,
with $r\in \Gamma_3(W)\otimes \Lambda^1$. The strategy used is to solve 
iteratively the equation 
$\delta(r)=R+\hat{\nabla}r -r^2/i\lambda$ with respect to the filtration
of $\Gamma(W)\otimes \Lambda$. It can be shown that the preceding equation
has a unique solution if the initial condition $\delta^{-1}(r)=0$ is satisfied.
It turns out that for such specific solution $r$ the Weyl curvature is nothing
else than $\Omega=-\omega$. Such a flat Abelian connection is called a 
\textsf{Fedosov connection} on $W$.

Given a Fedosov connection $D$ on $W$, let us see how it induces a star-product
on $M$.  The set of horizontal sections, i.e., sections $a$ such that $Da=0$, 
is an algebra for the product $\circ$. We shall denote this algebra 
by $\Gamma_D(W)$. The restriction of the projection 
$\sigma\colon\Gamma(W)\rightarrow C^\infty(M)[[\lambda]]$ to $\Gamma_D(W)$ is 
a bijective map. To any $f\in C^\infty(M)[[\lambda]]$ we can then associate 
a horizontal section denoted by $\sigma^{-1}(f)$. Then 
$f\star g= \sigma(\sigma^{-1}(f)\circ\sigma^{-1}(g))$ defines a star-product
on $M$.

The method of Fedosov extends to regular Poisson manifolds. An earlier
proof along the lines of De~Wilde and Lecomte and Omori, Maeda, and Yoshioka 
was published in 1992  by Masmoudi~\cite{Ma92} where the existence of 
tangential star-products is established.
\subsection{Classification, characteristic classes, and closed star-products}

\noindent\textbf{3.2.1 \ Classification.}
There are several occurrences in the literature of the second de Rham 
cohomology space $H^2_{dR}(M)$ in relation with the question of equivalence of 
star-products and deformed Poisson brackets. 
 
For formal 1-differentiable deformations of the Poisson bracket 
(i.e., formal Poisson bracket) on a symplectic manifold,
the first of these occurrences goes back to~\cite{BFFLS},
where it is shown that the obstructions to equivalence
(at each step in the power of the deformation parameter)
 are in $H^2_{dR}(M)$. By taking advantage of graded Lie structures related
to deformations, Lecomte~\cite{Le87} showed that equivalence
classes of formal Poisson brackets are in one-to-one correspondence with 
$H^2_{dR}(M)[[\lambda]]$.

There are similar occurrences for star-products. For example, it is implicit
in the work of Gutt~\cite{Gu79} and in \cite{GuIHP, Gu80} where the second 
and third Chevalley--Eilenberg cohomology spaces have been computed. Also,
as a consequence of their proof of existence when $b_3(M)=0$, 
Nero\-sla\-vsky  and  Vlasov~\cite{NV81} found that the classification for 
this class of star-products was given by sequences in $H^2_{dR}(M)$. 

The Fedosov construction provides a more geometric point of view 
on the classification problem. Fedosov's method allows to 
canonically construct a star-product on $(M,\omega)$
whenever a symplectic connection $\nabla$ is chosen. 
Any star-product obtained in such a way is  called a 
\textsf{Fedosov star-product}. The iterative construction reviewed in 
Section~\ref{fedosov} can be generalized to the case where one starts
with a formal symplectic form on $M$, 
$\omega_\lambda = \omega + \sum_{k\geq1} \lambda^k \omega_k$, the $\omega_k$'s 
being closed 2-forms on $M$. Hence we have a canonical method to construct 
a star-product $\star$ from the object $(M,\omega_\lambda,\nabla)$.
Then the Weyl curvature is equal to $-\omega_\lambda$ and Fedosov
has established that two Fedosov star-products $\star$ and $\star'$ 
are equivalent if and only if their Weyl curvatures belong to the same 
cohomology class in $H^2_{dR}(M)[[\lambda]]$. 

Shortly after Fedosov's paper, Nest and Tsygan~\cite{NT95} showed that
any differentiable star-product on a symplectic manifold is equivalent to
a Fedosov star-product. This fact allowed them to define the characteristic
class of a star-product~$\star$ to be the class of the Weyl curvature in
$H^2_{dR}(M)[[\lambda]]$ of any Fedosov star-product equivalent to~$\star$.
As a consequence, one gets a parametrization of the equivalence classes of
star-products on $(M,\omega)$ by elements in $H^2_{dR}(M)[[\lambda]]$.

The parametrization of equivalence classes of star-products has also been
made explicit by Bertelson, Cahen and Gutt~\cite{BCG97} (using older 
results in Gutt's thesis~\cite{Gu80}) and by Xu~\cite{Xu98}; the relation 
between classes of star-products and classes of 1-differentiable deformations 
of a Poisson bracket was described in \cite{Bon}. Deligne~\cite{De95} obtained
the same classification result, without any reference to a particular method 
of construction of deformations, by using cohomological \v Cech techniques 
and gerbes. Deligne's approach has been studied further by Gutt and 
Rawnsley~\cite{Gu00,GR99}. 
\vspace*{3mm}

\noindent\textbf{3.2.2 \ Trace and closed star products.}
The prominent role played by traces of operators  in index theory 
has its counterpart in deformation quantization. One can naturally define
a trace on the deformed algebra of smooth functions having compact support
$(C_c^\infty(M)[[\lambda]],\star)$.  It is  a 
$\mathbb{C}[[\lambda]]$-linear map $\mathrm{Tr}\colon (C_c^\infty(M)[[\lambda]] 
\rightarrow \mathbb{C}[\lambda^{-1},\lambda]]$
such that $\mathrm{Tr}(f\star g) = \mathrm{Tr}(g\star f)$. The negative powers 
of $\lambda$ appear already in the classical examples and $\mathrm{Tr}$ 
is allowed to take its values not only in formal series in $\lambda$ 
but in Laurent series in $\lambda^{-1}$ as well. 
Such a trace always exists and is essentially unique~\cite{NT95,Fe3}.
Around 1995, new algebraic versions of the (generalized) Atiyah-Singer 
theorem appeared: one is due to Nest and Tsygan~\cite{NT95} where 
Gelfand-Fuks techniques play an important role, the other, by Fedosov, 
uses Thom isomorphism. Recently, Fedosov~\cite{Fe00} published a new proof 
based on the topological invariants of symplectic connections 
of Tamarkin~\cite{Ta95}.
\vspace{1mm}

\noindent\textsl{3.2.2.1. \ Closed star-products.}
The trace allows to distinguish a particular class of star-products on a 
symplectic manifold $M$. These are the so-called \textsf{closed star-products}
introduced in~\cite{CFS}. A star-product $\star$ is closed if the map
$$
T\colon f \mapsto \frac{1}{(2\pi \lambda)^{m}}\int f_m \frac{\omega^m}{m!}
$$
is a trace ($2m$ is the dimension of $M$ and $f_m$ is the coefficient
of $\lambda^m$  in $f\in C_c^\infty(M)[[\lambda]]$). Closed star-products
exist on any symplectic manifold~\cite{OMY92} and any star-product is
equivalent to a closed star-product. Cyclic cohomology replaces the Hochschild
one for closed star-products: obstructions to existence and equivalence live
in the third and second  cyclic cohomology spaces, correspondingly. 
Another important feature of closed star-products is the existence of 
an invariant $\phi$ defined in~\cite{CFS}. It is a cocycle for the 
differential of the cyclic bicomplex of $M$, i.e., 
Hochschild in one direction and cyclic in the other. 
Let $\theta(f,g)= f\star g -f g$; for an integer $k$ define the map
$$
\phi_{2k}\colon C_c^\infty(M)[[\lambda]]^{\otimes 2k+1}
\rightarrow {\mathbb R},\quad (f_0,\ldots,f_{2k})\mapsto
\tau(f_0\star \theta(f_1,f_2)\star\cdots\star
\theta(f_{2k-1},f_{2k}))
$$
where $\tau=(2\pi\lambda)^{m}m! T$ does not depend on $\lambda$. It turns out
that $\phi_{2k}=0$ for $k>m$ and the cocycle is defined by
$\phi=\sum_k \phi_{2k}/k!$. The computation of $\phi$, called the 
CFS-invariant, has been performed for the case where $M$ is a cotangent 
bundle of a compact Riemannian manifold~\cite{CFS}. Notice that in that 
case, associative deformations correspond to the standard calculus of 
pseudo-differential operators, and the CFS-invariant is nothing else than 
the Todd class of $M$. This result has been generalized by 
Halbout~\cite{Hal99} to the case of any symplectic manifold.

\section{Metamorphoses}
\addtocounter{equation}{-2}
\subsection{Poisson manifolds and Kontsevich's formality} 

The problem of deformation quantization of general Poisson differentiable 
manifolds has been left open for a long time. 
For the nonregular Poisson case, first examples of star-products appeared
in \cite{BFFLS} in relation with the quantization of angular 
momentum. They were defined on the dual of $\mathfrak{so}(n)$ endowed with its
natural Kirillov--Poisson structure. The case for any Lie algebra follows from 
the construction given by Gutt \cite{Gu83}
of a star-product on the cotangent
bundle of a Lie group~$G$. This star-product restricts to a star-product on
the dual of the Lie algebra of~$G$. It translates the associative structure
of the enveloping algebra in terms of functions on the dual of the Lie
algebra of~$G$. Omori, Maeda and Yoshika~\cite{OMY94} have constructed 
quantization for a class of quadratic Poisson structures.

There is no equivalent of Darboux theorem for Poisson structures and this fact
is at the origin of the main difficulty to construct (even local) star-products
in the Poisson case. The situation changed drastically in 1997, when 
Kontsevich~\cite{Konsi97b} proved his Formality Conjecture~\cite{Konsi97a}.
The existence of star-products on any smooth Poisson manifold appears as a 
consequence of the Formality Theorem.  

We shall review the construction of a star-product on $\mathbb{R}^d$ given in 
\cite{Konsi97b} and then present the  essential steps of the Formality Theorem.
\vspace*{3mm}

\noindent\textbf{4.1.1 \ Kontsevich star-product.} 
Consider $\mathbb{R}^d$ endowed with a Poisson bracket $\alpha$. 
We denote by  $(x^1,\ldots,x^d)$ the coordinate system on $\mathbb{R}^d$;
the Poisson bracket of two smooth functions  $f,g$ is given by 
$\alpha(f,g)= \sum_{1\leq i,j\leq n}\alpha^{ij}\partial_i f \partial_j g$, 
where $\partial_k$ denotes the 
partial derivative with respect to $x^k$. 
Instead of the whole of  $\mathbb{R}^d$, one may consider an open subset of it.

The formula for the Kontsevich star-product is conveniently defined
by considering, for each $n\geq 0$, a family of oriented graphs $G_{n}$. To a 
graph $\Gamma \in G_{n}$ is associated a bidifferential operator $B_\Gamma$ and
a weight $w(\Gamma)\in \mathbb{R}$. 
The sum $\sum_{\Gamma\in G_{n}}
w(\Gamma) B_\Gamma$ gives  the coefficient of $\lambda^n$, i.e., the cochain 
$C_n$ of the star-product. 

For later use, we shall reproduce the general definitions of~\cite{Konsi97b}.
They are more general than the ones needed in this section for the star-product.

An oriented graph $\Gamma$ belongs to $G_{n,m}$, 
with $n,m\geq 0$, $2n-m+2\geq0$, if:
\begin{itemize}
\item[i)] The set of vertices $V_\Gamma$ of $\Gamma$ has $n+m$ elements 
labeled $\{1,\ldots,n;\bar{1},\ldots,\bar{m}\}$.
The vertices $\{1,\ldots,n\}$ are said of the first type and 
$\{\bar{1},\ldots,\bar{m}\}$ are said of the second type.

\item[ii)]  $E_\Gamma$ is the set of oriented edges of $\Gamma$. 
$E_\Gamma$ has $2n-m+2$ edges. There is
no edge starting at a vertex of the second type. The set of edges starting 
at a vertex of the first type $k$ is denoted by $\mathrm{Star}(k)$ and its 
cardinality by $\sharp(k)$, hence 
$\sum_{1\leq k\leq n}  \sharp(k)=2n+m-2$.
The edges of $\Gamma$ starting at vertex $k$ will be denoted by 
$\{e_k^1,\ldots,e_k^{\sharp(k)}\}$. When it is needed to make explicit the  
starting and ending vertices of some edge $v$, we shall write $v=(s(v),e(v))$, 
where $s(v)\in \{1,\ldots,n\}$ is its starting vertex and
$e(v) \in \{1,\ldots,n;\bar{1},\ldots,\bar{m}\}$ is its ending vertex.

\item[iii)] $\Gamma$ has no loop (edge starting at some vertex and ending at
that vertex) and no parallel multiple edges (edges sharing the same starting
and ending vertices).
\end{itemize}

Let us specialize the preceding definition to the case of star-products. Here 
$m=2$ and only the family of graphs $G_{n,2}$, $n\geq0$, will be of interest.
Moreover there are exactly $2$ edges starting at each vertex of the first type,
i.e., $\sharp{k}=2$, $\forall k \in \{1,\ldots,n\}$. Hence $E_\Gamma$ has $2n$ 
edges for $\Gamma\in G_{n,2}$. 
For notational convenience we shall write the pair of edges 
$\{e^{1}_k,e^{2}_k\}$ as $\{i_k,j_k\}$. Also $\{i_k,j_k\}$ will
stand for indices of the Poisson tensor $\alpha^{i_kj_k}$, 
with $1\leq i_k, j_k\leq d$.
The set of graphs in $G_{n,2}$ has $(n(n+1))^n$ elements. 
For $n=0$, $G_{0,2}$ has only one element: 
the graph having as set of vertices $\{\bar1,\bar2\}$ and no edges.

A bidifferential operator $(f,g)\mapsto B_\Gamma(f,g)$,
$f,g\in C^\infty(\mathbb{R}^d)$,
is associated to each graph $\Gamma\in G_{n,2}$, $n\geq1$, in the following way. 
To each vertex~of the first type $k$ from where the edges $\{i_k,j_k\}$
start, one associates the components
$\alpha^{i_kj_k}$ of the Poisson tensor; $f$ is associated to the vertex~$\bar1$ 
and $g$ to the vertex $\bar2$. Each edge, e.g. $i_k$, acts by partial 
differentiation with respect to $x^{i_k}$ on its ending vertex. 

For example, the graph in Fig.~1 gives the bidifferential operator
\begin{figure}\label{fig1}
\centerline{\epsfig{file=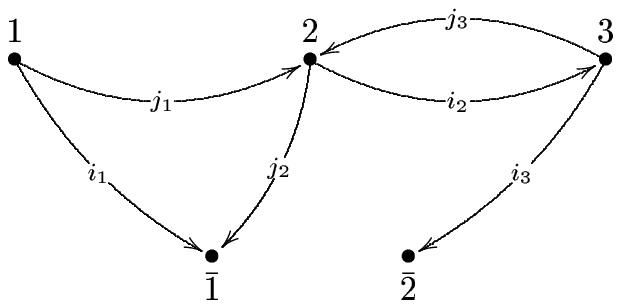,width=60mm}}
\caption{A graph in $G_{3,2}$.}
\end{figure}
$$
B_\Gamma(f,g) = \sum_{0\leq i_\ast,j_\ast\leq d}
\alpha^{i_1j_1}\,\partial_{j_1j_3}\alpha^{i_2j_2}\,
\partial_{i_2}\alpha^{i_3j_3}\,\partial_{i_1j_2}f\,\partial_{i_3}g.
$$
Notice that for $n=0$, we simply have the usual product of $f$ and $g$.
The general formula for $B_\Gamma(f,g)$ for $\Gamma\in G_{n,2}$ is given by
a sum over all maps $I\colon E_\Gamma\rightarrow \{1,\ldots,d\}$, 
\begin{equation}\label{polydif}
B_\Gamma(f,g) = \sum_{I} \left(
\left(\prod_{k=1}^n\prod_{k'=1}^n 
\partial_{I(k',k)} \alpha^{I(e^1_k) I(e^2_k)}\right)
\left(\prod_{k_1=1}^n \partial_{I(k_1,\bar1)} f\right)
\left(\prod_{k_2=1}^n \partial_{I(k_2,\bar2)} g\right)
\right),
\end{equation}
where the usual convention for missing factors is used (i.e., whenever
an edge is absent in~$E_\Gamma$, the corresponding factor is 1).

Finally the weight $w(\Gamma)$ of a graph $\Gamma$ is defined as follows
(a general definition will be given in Section~4.1.2.
Let  $\mathcal{H}=\{z\in\mathbb{C}\,|\,\mathrm{Im}(z)>0\}$ 
be the upper half-plane with its Lobachevsky hyperbolic metric.
The configuration space of $n$ distinct points in $\mathcal{H}$ is denoted
by $\mathcal{H}_n$. It is an open submanifold of $\mathbb{C}^n$ with its
natural orientation. 

Let $\phi\colon \mathcal{H}_2\rightarrow \mathbb{R}/2\pi \mathbb{Z}$ be the 
function (angle function):
\begin{equation}\label{phih}
\phi(z_1,z_2)=\frac{1}{2\sqrt{-1}}\mathrm{Log}
\Bigl(\frac{(z_2-z_1)(\bar{z}_2-z_1)}
{(z_2-\bar{z}_1)(\bar{z}_2-\bar{z}_1)}\Bigr).
\end{equation}
$\phi(z_1,z_2)$ is extended by continuity for $z_1,z_2\in \mathbb{R}$, 
$z_1\neq z_2$. The function $\phi$ measures the angle at $z_1$ between the 
geodesic passing by $z_1$ and the point at infinity, and the geodesic 
passing by $z_1$ and $z_2$ for the hyperbolic metric.

For a graph $\Gamma\in G_{n,2}$, the vertex $k$, $1\leq k\leq n$, is 
associated to a variable  $z_k\in\mathcal{H}$, the vertex $\bar1$
to $0\in\mathbb{R}$, and the vertex $\bar2$ to $1\in\mathbb{R}$. Doing so,
we can associate a function $\tilde\phi_v$ on $\mathcal{H}_n$ to an edge 
$v\in E_\Gamma$ by $\tilde\phi_v = \phi(s(v),e(v))$.

The weight  $w(\Gamma)$  is defined by integrating a $2n$-form over 
$\mathcal{H}_n$ (in fact on some compactified space, 
cf. Section~4.1.2:
\begin{equation}\label{w}
w(\Gamma)=\frac{1}{n!(2\pi)^{2n}}\int_{\mathcal{H}_n}\bigwedge_{1\leq k\leq n}
\big(d\tilde\phi_{e_k^1}\wedge d\tilde\phi_{e_k^2}\big),
\end{equation}
where $\{e_k^1,e_k^2\}$, $1\leq k \leq n$, is the set of edges of $\Gamma$.
The weights are universal in the sense that they do not depend on the 
Poisson structure or the dimension~$d$. 

The above construction goes over to any formal Poisson bracket
$\alpha_\lambda = \sum_{k\geq0} \lambda^k \alpha_k$.

\begin{theorem}[Kontsevich~\cite{Konsi97b}]\label{thkonsi}
The map $\star\colon C^\infty(\mathbb{R}^d)\times C^\infty(\mathbb{R}^d)
\rightarrow C^\infty(\mathbb{R}^d)[[\lambda]]$
defined by
$$(f,g)\mapsto f\star g= \sum_{n\geq0} \lambda^n \sum_{\Gamma\in G_{n,2}}
w(\Gamma)B_\Gamma(f,g)
$$ 
defines a differential star-product on $(\mathbb{R}^d,\alpha)$.
The classes of equivalence of star-products are in a one-to-one correspondence
with the classes of equivalence of formal Poisson brackets.
\end{theorem}

\begin{figure}\label{fig2}
\centerline{\epsfig{file=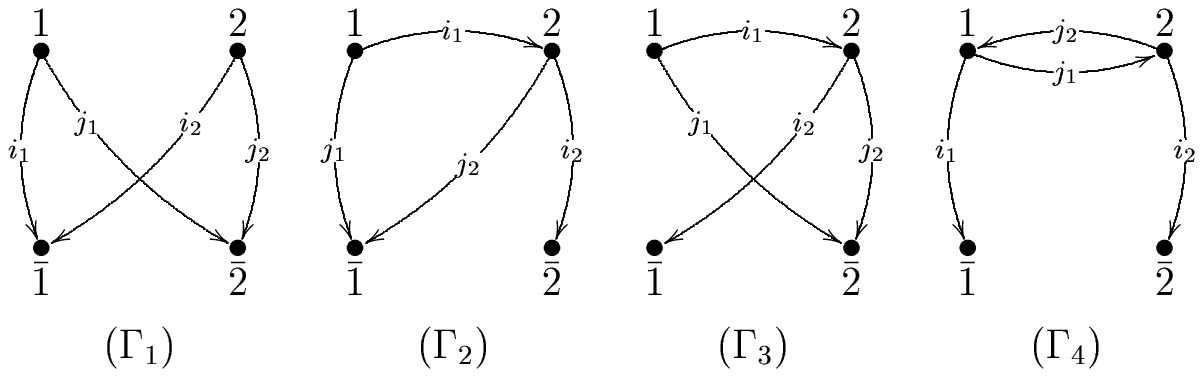,width=110mm}}
\caption{Graphs in $G_{2,2}$ contributing to $C_2$.}
\end{figure}

The first cochains of the Kontsevich star-product can be computed by elementary 
means. Of course we have $C_0(f,g)=fg$. At first order, $G_{1,2}$ has 2 graphs:
one has weight $\frac{1}{2}$ and as bidifferential operator the Poisson bracket
$\{f,g\}_\alpha = \alpha(df\wedge dg)$, the other graph has weight 
$-\frac{1}{2}$ and bidifferential operator $\{g,f\}_\alpha$. 
The sum gives $C_1(f,g)= \{f,g\}_\alpha$. 
All the graphs contributing to $C_2$, up to symmetries, are shown in~Fig.~2, 
and the second order term is given by:
\begin{eqnarray}\label{c2}
C_2(f,g)\hspace*{2mm} = &&\hspace*{-5mm}\frac{1}{2}\,\alpha^{i_1j_1}\,
\alpha^{i_2j_2}\,\partial_{i_1 i_2} f\, \partial_{j_1 j_2} g\nonumber\\
&&+ \frac{1}{3}\,\alpha^{i_1j_1}\,\partial_{i_1}\alpha^{i_2j_2}\,
(\partial_{j_1 j_2} f\,\partial_{i_2} g +\partial_{i_2} f\, 
\partial_{j_1 j_2}g)\\
&&-\frac{1}{6}\,\partial_{j_2}\alpha^{i_1j_1}\,
\partial_{j_1}\alpha^{i_2j_2}\,\partial_{i_1} f\, \partial_{i_2}g,\nonumber
\end{eqnarray}
where summation over repeated indices is understood. For higher cochains,
there is no systematic way to compute weights of graphs, especially
in the presence of cycles (e.g. graph $\Gamma_4$ in Fig.~2), but we should
mention that Polyak~\cite{Po01} has provided an interpretation of weights 
in term of degree of maps which allowed him to ``easily'' compute a large class
of graphs and all of them for the linear Poisson structure case.

The Kontsevich star-product has been extensively studied for linear Poisson 
structures, i.e., on the dual of a Lie algebra. An important
result in that context is due to Shoikhet~\cite{Sho01} who showed the
vanishing of the weight of all even wheel graphs. As a consequence, on the dual
of a Lie algebra, the Kontsevich star-product coincides with the one given
by the Duflo isomorphism. Moreover, the Kontsevich star-product induces
a convolution on the space of distributions supported at $0$ in the
Lie algebra $\mathfrak{g}$.  This has allowed \cite{ADS} to prove the  
conjecture on invariant distributions on a Lie group formulated by 
Kashiwara and Vergne in 1978. 

The linear Poisson case allows  to quantize the class of
quadratic Poisson brackets that are in the image of a map
introduced by Drinfeld which associates a quadratic Poisson 
to a linear Poisson bracket; not all quadratic brackets arise in this 
way~\cite{MMR}.
\vspace*{3mm}

\noindent\textbf{4.1.2 \ Formality Theorem.} 
Theorem~\ref{thkonsi} is a particular consequence of a much more general
result called Formality Theorem. Given a smooth  manifold $M$ one may
consider two algebraic structures constructed on
the associative algebra $A\equiv C^\infty(M)$, namely,
the Hochschild complex $C^\bullet(A,A)$ and its cohomology $H^\bullet(A,A)$.
They both have a structure of differential graded Lie algebra. 
Roughly speaking, the Formality Theorem states that these two differential 
graded Lie algebras are ``equivalent'' in a sense we now make precise.
\vspace{2mm}

\noindent\textsl{4.1.2.1. \ Differential graded Lie algebras.}~Let 
$\mathbb{K}$ be a field 
of characteristic zero. All vector spaces are over $\mathbb{K}$.
Let $V$ be a $\mathbb{Z}$-graded vector space:
$V=\oplus_{k\in\mathbb{Z}} V^k$. An element $x\in V$ is called homogeneous
of degree $k$ if $x\in V^k$. The degree of a homogeneous element
$x$ is denoted~$|x|$. From now on, ``graded'' will 
always mean ``$\mathbb{Z}$-graded''.

A \textsf{differential graded Lie algebra} over $\mathbb{K}$ 
(\textsf{DGLA} for short) is a graded vector space 
$\mathfrak{g} = \oplus_{k\in\mathbb{Z}} \mathfrak{g}^k$
endowed with Lie superbracket $[\cdot,\cdot]$ and a differential
$d$ of degree~$1$. It means that $(\mathfrak{g}, [\cdot,\cdot],d)$
fulfills ($a,b,c$ are homogeneous elements):
\begin{itemize}
\item[1)] $[\cdot,\cdot]$ is a bilinear map, 
$[\cdot,\cdot]\colon \mathfrak{g}^i\times 
\mathfrak{g}^j\rightarrow \mathfrak{g}^{i+j}$,
satisfying:

$[a,b] = - (-1)^{|a| |b|} [b,a]$,

$(-1)^{|a||c|}[a,[b,c]] + (-1)^{|b||a|}[b,[c,a]] + (-1)^{|c||b|}[c,[a,b]] =0$.

\item[2)] $d\colon \mathfrak{g}^i\rightarrow \mathfrak{g}^{i+1}$ 
is a linear map satisfying $d^2=0$ and:

$ d[a,b] = [da,b] + (-1)^{|a|}[a,db]$.
\end{itemize}

The Hochschild complex $C^\bullet(A,A)=\oplus_{k\geq0} C^k(A,A)$ 
consists of polydifferential operators on $M$, 
$C^k(A,A)=\{D\, |\, D\colon A^{\otimes^k}\rightarrow A\}$. 
In the following, we will detail the ingredients that make the
Hochschild complex a \textsf{DGLA}. We shall look at the Hochschild 
complex as a  graded vector space 
$D^\bullet_\mathrm{poly}(M)=\oplus_{k\in\mathbb{Z}}D^k_\mathrm{poly}(M)$, 
where $D^k_\mathrm{poly}(M)$ is equal to $C^{k+1}(A,A)$ for $k\geq -1$ 
and equal to $\{0\}$ otherwise. Notice the shift by one unit in the grading
of $D^\bullet_\mathrm{poly}(M)$ with respect to that of $C^\bullet(A,A)$:
$|D|=k$ if $D\in C^{k+1}(A,A)$.

The \textsf{Gerstenhaber bracket} $[\cdot,\cdot]_G$ on 
$D^\bullet_\mathrm{poly}(M)$ is defined for 
$D_i \in D^{k_i}_\mathrm{poly}(M)$ by:
$$
[D_1, D_2]_G = D_1\circ D_2 - (-1)^{k_1 k_2}D_2\circ D_1,
$$
where $\circ\colon D^{k_1}_\mathrm{poly}(M)\times D^{k_2}_\mathrm{poly}(M)
\rightarrow D^{k_1+k_2}_\mathrm{poly}(M)$ 
is a nonassociative composition law for polydifferential operators:
\begin{equation*}
\begin{split}
(D_1\circ D_2)&(f_0,\ldots, f_{k_1+k_2}) \\
&=
\sum_{0\leq j\leq k_1} (-1)^{jk_2} D_1(f_0,\ldots, f_{j-1}, D_2(f_j,
\ldots,f_{j+k_2}), f_{j+k_2+1},\ldots f_{k_1+k_2}).
\end{split}
\end{equation*}
Finally, we have the usual Hochschild differential $b$ defined by~(\ref{bH}).
It can be expressed in terms of the Gerstenhaber bracket. 
If $m\colon A\otimes A\rightarrow A$ denotes
the product of $A$ (i.e., usual product of functions on $M$), then one can
verify that $b = - [\ \cdot\ ,m]_G$. The usual  Hochschild differential is not
a differential of the Gerstenhaber bracket, but $\delta = [m,\ \cdot\ ]_G$ is
a differential of degree $1$ (for $D\in D^{k}_\mathrm{poly}(M)$, 
we have $bD=(-1)^k \delta D$).
It is an important fact that 
$(D^\bullet_\mathrm{poly}(M), [\cdot,\cdot]_G, \delta)$ is a \textsf{DGLA}.

The cohomology $H^\bullet(A,A)$ of the complex $C^\bullet(A,A)$ 
(with $b$ or $\delta$ as differential) is the space of  polyvectors on $M$:
$\Gamma(\wedge^\bullet TM)$~\cite{Ve, HKR}. We denote by 
$T^\bullet_\mathrm{poly}(M)$ the graded vector space 
$T^\bullet_\mathrm{poly}(M)=\oplus_{k\in\mathbb{Z}}T^k_\mathrm{poly}(M)$, 
where $T^k_\mathrm{poly}(M) = \Gamma(\wedge^{k+1} TM)$ for $k\geq -1$ 
and is equal to $\{0\}$ otherwise. On $T^\bullet_\mathrm{poly}(M)$
one has the Schouten--Nijenhuis bracket $[\cdot,\cdot]_{SN}$ which is 
the natural extension of the bracket of vector fields.
On decomposable tensors it is defined by:
\begin{equation}\label{schouten}
\begin{split}
[&\xi_0\wedge\cdots\wedge\xi_k,\eta_0\wedge\cdots\wedge\eta_l]_{SN}\\
&=\sum_{0\leq i\leq k}\sum_{0\leq j\leq l} (-1)^{i+j} [\xi_i,\eta_j]
\wedge\xi_0 \wedge\cdots\wedge\hat\xi_i\wedge\cdots\wedge\xi_k\wedge
\eta_0\wedge\cdots\wedge\hat\eta_j\wedge\cdots\wedge\eta_l.
\end{split}
\end{equation}
The Schouten--Nijenhuis bracket makes $(T^\bullet_\mathrm{poly}(M), 
[\cdot,\cdot]_{SN})$ into a graded Lie algebra. On it, take as the  
differential $0$, then $(T^\bullet_\mathrm{poly}(M), [\cdot,\cdot]_{SN},0)$ 
is trivially a \textsf{DGLA}.
\vspace{2mm}

\noindent\textsl{4.1.2.2. $L_\infty$-morphism.}
We recall some  definitions needed to state the Formality Theorem. 
A \textsf{graded coalgebra} $\mathfrak{c}$ is a graded
vector space  with a coproduct $\Delta\colon \mathfrak{c}_i \rightarrow
\oplus_{m+n=i} \mathfrak{c}_m\otimes\mathfrak{c}_n$. To any graded vector
space $\mathfrak{g}$ is associated a graded coalgebra 
$\mathcal{C}(\mathfrak{g}) = \otimes_{k\geq1} (\wedge^k\mathfrak{g})[k]$,
where $[k]$ indicates the shift to the left by $k$ units in the grading.

An \textsf{$L_\infty$-algebra} $(\mathfrak{g},Q)$ is a graded vector space 
$\mathfrak{g}$ with a codifferential $Q$ satisfying $Q^2=0$
on $\mathcal{C}(\mathfrak{g})$.
$L_\infty$-algebras are also known as \textsf{strong(ly) homotopy Lie 
algebras} (see e.g. \cite{Sta}). 
The differential $Q$ defines a sequence $\{Q_k\}$ of morphisms
$Q_k\colon \wedge^k(\mathfrak{g}) \rightarrow \mathfrak{g}[2-k]$, $k\geq1$
and the relation $Q^2=0$ puts quadratic constraints on the $Q_k$'s. These 
constraints  imply that $Q_1\colon \mathfrak{g}\rightarrow \mathfrak{g}[1]$ 
is a differential of degree 1, hence $(\mathfrak{g},Q_1)$ is a complex, and 
$Q_2\colon \wedge^2(\mathfrak{g}) \rightarrow \mathfrak{g}$
defines a superbacket on $\mathfrak{g}$ for which $Q_1$ is compatible.
Hence any $L_\infty$-algebra $(\mathfrak{g},Q)$ such that $Q_k=0$ for $k\geq3$
is a~\textsf{DGLA}.

An $L_\infty$-morphism $\mathcal{F}$ between two $L_\infty$-algebras
$(\mathfrak{g}_1,Q^{(1)})$ and $(\mathfrak{g}_2,Q^{(2)})$ is a morphism
between the associated coalgebras such that 
$\mathcal{F}Q^{(1)} = Q^{(2)}\mathcal{F}$.
We have similar constraints on the sequence of morphisms $\{\mathcal{F}_k\}$:
$
\mathcal{F}_k\colon \wedge^k(\mathfrak{g_1}) 
\rightarrow \mathfrak{g_2}[1-k],\quad k\geq1.
$
At first order, they imply that the formality $\mathcal{F}$ is a morphism of 
complexes $(\mathfrak{g}_1,Q^{(1)}_1) \rightarrow (\mathfrak{g}_2,Q^{(2)}_1)$,
but the constraint at order 2 shows that $\mathcal{F}_2$ is not a homomorphism
of \textsf{DGLA}.
Finally, an $L_\infty$-morphism is a quasi-morphism
if the morphism of complexes
$\mathcal{F}_1\colon (\mathfrak{g}_1,Q^{(1)}_1) 
\rightarrow (\mathfrak{g}_2,Q^{(2)}_1)$ induces an isomorphism
between the corresponding cohomology spaces (i.e., quasi-isomorphism
of complexes in the usual sense).

The two \textsf{DGLA} $(T^\bullet_\mathrm{poly}(M), [\cdot,\cdot]_{SN},0)$ and 
$(D^\bullet_\mathrm{poly}(M), [\cdot,\cdot]_G, \delta)$ are $L_\infty$-algebras
by definition. The natural injection 
$\mathcal{I}\colon T^\bullet_\mathrm{poly}(M)
\rightarrow D^\bullet_\mathrm{poly}(M)$ defined by
$$
\mathcal{I}(\xi_1\wedge\cdots\wedge\xi_k).(f_1,\ldots,f_k)
= \frac{1}{k!} \det(\xi_i(f_j)),
$$
is a quasi-isomorphism of complexes, but is not a homomorphism
of \textsf{DGLA} since it is incompatible with the brackets.
But Kontsevich has shown that:
\begin{theorem}[Formality]
There is a quasi-isomorphism 
$$
\mathcal{K}\colon(T^\bullet_\mathrm{poly}(M), [\cdot,\cdot]_{SN},0)\rightarrow
(D^\bullet_\mathrm{poly}(M), [\cdot,\cdot]_G, \delta)
$$ 
such that $\mathcal{K}_1$ is equal to $\mathcal{I}$.
\end{theorem}
The main part of the proof consists in writing down  
explicitly a  quasi-isomorphism  for the case $M=\mathbb{R}^d$.
The quasi-isomorphism $\mathcal{K}$ is expressed in terms of graphs and weights 
and it generalizes the formula for the star-product given in Section~4.1.1.
\vspace{2mm}

\noindent\textsl{4.1.2.3.\ Construction of} $\mathcal{K}$ \textsl{for} 
$M=\mathbb{R}^d$.
We have already defined the family of graphs $G_{n,m}$. Here we need
a generalization of the bidifferential operators $B_\Gamma$ of the star-product
formula. For a graph $\Gamma$ in $G_{n,m}$, take $(\alpha_1,\ldots, \alpha_n)$
$n$ polyvectors such that $\alpha_k\in T^{\sharp(k)-1}_\mathrm{poly}(M)$. 
By definition of $G_{n,m}$ we have $\sum_{1\leq k\leq n}  \sharp(k)=2n+m-2$,
and consequently $\alpha_1\wedge\cdots\wedge\alpha_n$ is in 
$T^{n+m-2}_\mathrm{poly}(M)$. Define a polydifferential
operator $B_\Gamma(\alpha_1\wedge\cdots\wedge\alpha_n)$ belonging
to $T^{m-1}_\mathrm{poly}(M)$ by a formula similar to~(\ref{polydif}). For 
$m$ smooth functions $(f_1,\ldots, f_m)$ on $\mathbb{R}^d$, we have:
\begin{equation}\label{polydiff}
\begin{split}
B_\Gamma(\alpha_1\wedge&\cdots\wedge\alpha_n).
(f_1,\ldots, f_m) = \\
&\sum_{I} \left(
\left(\prod_{k=1}^n\prod_{k'=1}^n 
\partial_{I(k',k)} \alpha_k^{I(e^1_k)\cdots I(e^{\sharp(k)}_k)}\right)
\left(\prod_{i=1}^m\prod_{j=1}^n \partial_{I(j,\bar i)} f_i\right)
\right),
\end{split}
\end{equation}
where the sum is on all maps $I\colon E_\Gamma\rightarrow \{1,\ldots,d\}$.

The weight of a graph  is defined by integration on 
compactified configuration spaces that we now make precise.
For $2n+m-2\geq0$, let
$$\mathsf{Conf}_{n,m}=\{(p_1,\ldots,p_n;q_1,\ldots,q_m)\in 
\mathcal{H}^n\times\mathbb{R}^m\ |\ 
p_i\neq p_j\ \mathrm{for} \ i\neq j; q_i \neq q_j\ \mathrm{for} \ i\neq j\}.
$$
$\mathsf{Conf}^+_{n,m}$ will denote the connected
component of $\mathsf{Conf}_{n,m}$ for which $q_1<\ldots<q_m$.
$\mathsf{Conf}_{n,m}$ is acted upon by the group $G$ of transformations
$z\mapsto az+b$, $a>0$, $b\in\mathbb{R}$. 
Define $\mathsf{C}_{n,m}=\mathsf{Conf}_{n,m}/G$: It inherits
a natural orientation. $\mathsf{C}^+_{n,m}$ will denote the connected
component corresponding to $\mathsf{Conf}^+_{n,m}$. The compactified spaces
$\bar{\mathsf{C}}^+_{n,m}$ are manifolds with corners. 
The weight of a graph $\Gamma\in G_{n,m}$ is defined by:
\begin{equation}\label{wa}
w(\Gamma)=\prod_{1\leq k \leq n}\frac{1}{(\sharp{k})!}\frac{1}{(2\pi)^{2n+m-2}}
\int_{\bar{\mathsf{C}}^+_{n,m}}
\bigwedge_{e\in E_\Gamma} d\tilde\phi_{e}\big.
\end{equation}
Notice that the order in the set of edges $E_\Gamma$ gives the order 
in the wedge product. We can now write down the sequence of morphisms 
$\{ \mathcal{K}_n\}$ of the quasi-isomorphism $\mathcal{K}$:
$$
\mathcal{K}_n =\sum_{m\geq0}\sum_{\Gamma\in G_{n,m}}w(\Gamma)B_\Gamma,
\quad n\geq1.
$$
The verification that $\mathcal{K}$ is a quasi-isomorphism amounts to check 
an infinite number of quadratic equations for the weights, which are 
essentially a consequence of Stokes theorem applied to the weights. 
One has to keep track of all the signs mixed by the gradings and the 
orientations of the configuration spaces. A detailed description of these 
computations can be found in~\cite{AMM00}.

\subsection{Field theory and path integral formulas} 
\vspace*{3mm}

\noindent\textbf{4.2.1 \ Star-products in infinite dimension.}
The deformation quantization of a given classical field theory consists in
giving a proper definition for a star-product on the infinite-dimensional
manifold of initial data for the classical field equation  and
constructing with it, as rigorously as possible, whatever physical
expressions are needed.
As in other approaches to field theory, here also one faces serious divergence
difficulties as soon as one is considering interacting fields theory, and
even at the free field level if one wants a mathematically rigorous theory.
But the philosophy in dealing with the divergences is significantly different
and one is in position to take advantage of the cohomological features of
deformation theory to perform what can be called {\it cohomological
renormalization}~\cite{Di93}.

Poisson structures are known on infinite-dimensional manifolds since
a long time and there is an extensive literature on this subject. 
A typical structure, for our purpose, is a weak
symplectic structure such as that defined by Segal \cite{Se74} 
(see also \cite{Ko74})
on the space of solutions of a classical field equation like
$\Box\Phi=F(\Phi)$, where $\Box$ is the d'Alembertian. Now if one considers
scalar-valued functionals $\Psi$ over such a space of solutions, i.e., over
the phase space of initial conditions $\varphi(x)=\Phi(x,0)$ and
$\pi(x)={\frac{\partial}{\partial t}}\Phi(x,0)$, one can consider a Poisson
bracket defined by
\begin{equation} \label{Poissonf}
P(\Psi_1,\Psi_2)=\int\left({\frac{\delta\Psi_1}{\delta\varphi}}
{\frac{\delta\Psi_2}{\delta\pi}}-{\frac{\delta\Psi_1}{\delta\pi}}
{\frac{\delta\Psi_2}{\delta\varphi}}\right)dx
\end{equation}
where $\delta$ denotes the functional derivative. The problem is that while
it is possible to give a precise mathematical meaning to (\ref{Poissonf}) by
specifying an appropriate algebra of functionals, the
formal extension to powers of $P$, needed to define the Moyal bracket,
is highly divergent, already for $P^2$.
This should not be so surprising for physicists who know from experience that
the correct approach to field theory starts with  normal ordering, and that
there are infinitely many inequivalent representations of the canonical
commutation relations (as opposed to the von Neumann uniqueness in the
finite-dimensional case, for projective representations), even if in recent
physical literature some are working formally with Moyal product.

Starting with some star-product $\star$ (e.g., an infinite-dimensional version
of a Moyal-type product or, better, a star-product similar to the
normal star-product (\ref{np})) on the manifold of initial data, one would
interpret various divergences appearing in the theory in terms of
coboundaries (or cocycles) for the relevant Hochschild cohomology.
Suppose that we are suspecting that a term in a cochain of the product $\star$
is responsible for the appearance of divergences. Applying an iterative 
procedure of equivalence, we can try to eliminate it, or at least get a 
lesser divergence, by subtracting at the relevant order a divergent 
coboundary; we would then get a better theory with a new star-product, 
``equivalent'' to the original one.
Furthermore, since in this case we can expect to have at each order an
infinity of non equivalent star-products, we can try to subtract a cocycle
and then pass to a nonequivalent star-product whose lower order cochains
are identical to those of the original one. We would then make an analysis of
the divergences up to order $\hbar^r$, identify a divergent cocycle, remove
it, and continue the procedure (at the same or hopefully a higher order).
Along the way one should preserve the usual properties of a quantum field
theory (Poincar\'e covariance, locality, etc.) and the construction of
adapted star-products should be done accordingly. The complete implementation
of this program should lead to a cohomological approach to renormalization
theory. It would be interesting to formulate the Connes--Kreimer~\cite{CoKr} 
rigorous renormalization procedure in a way that would fit in this pattern.

A very good test for this approach would be to start from classical
electrodynamics, where (among others) the existence of global
solutions and a study of infrared divergencies were recently rigorously
performed \cite{FST97}, and go towards mathematically rigorous QED.
Physicists may think that spending so much effort in trying to give complete 
mathematical sense to recipes that work so well is a waste of time, but the
fact it has for so long resisted a rigorous formulation is a good indication
that the new mathematical tools needed will prove very efficient.
\vspace{2mm}

\noindent\textsl{4.2.1.1.\ Normal star-product.}~In the case of free fields, 
one can write down an explicit expression for
a star-product corresponding to normal ordering. Consider a (classical)
free massive scalar field $\Phi$ with initial data $(\phi,\pi)$ in the
Schwartz space ${\cal{S}}$. The initial data  $(\phi,\pi)$ can
advantageously be replaced by their Fourier modes $({\bar{a}},a)$, which
for a massive theory are also in ${\cal{S}}$ seen as a real vector space.
After quantization $({\bar{a}},a)$ become 
the usual creation and annihilation operators, respectively. 
The normal star-product $\star_N$ is formally equivalent
to the Moyal product and an integral 
representation for $\star_N$ is given by~\cite{Di90}:
\begin{equation} \label{np}
(F\star_N G)(\bar{a},a)=\int_{{\cal{S}}'\oplus {\cal{S}}'}
d\mu(\bar{\xi},\xi) F(\bar{a},a +\xi) G(\bar{a} + \bar{\xi}, a),
\end{equation}
where, by Bochner--Minlos theorem, $\mu$ is the Gaussian measure on 
${\cal{S}}'\oplus {\cal{S}}'$ defined by the characteristic function 
$\exp(-\frac{1}{\hbar}\int dk\ {\bar{a}}(k)a(k))$
and $F,G$ are holomorphic functions with semi-regular kernels. 
The reader may object that creation and annihilation operators are 
operator-valued distributions and one should consider $({\bar{a}},a)$ in
$\mathcal{S}'\oplus \mathcal{S}'$ instead. In fact the distribution aspect is
present in the very definition of the cochains of the 
star-product through (\ref{np}).

Likewise, Fermionic fields can be cast in that framework by considering 
functions valued in some infinite dimensional Grassmann algebra and 
super-Poisson brackets (for the deformation quantization of the latter 
see e.g. \cite{Bm96}).

For the normal product (\ref{np}) one can formally consider interacting fields.
It turns out that the star-exponential of the Hamiltonian is, up to a
multiplicative well-defined function, equal to Feynman's path integral.
For free fields, we have a mathematical meaningful equality between
the star-exponential and the path integrals as both of them are defined
by a Gaussian measure, and hence well-defined. In the interacting fields
case, giving a rigorous meaning to either of them would give a meaning
to the other.

The interested reader will find in~\cite{Di90} calculations performing some
steps in the above direction, for free scalar fields and the Klein--Gordon
equation, and an example of cancellation of some infinities in
$\lambda\phi_2^4$-theory via a $\lambda$-dependent star-product formally
equivalent to a normal star-product. 

Recently many field theorists became interested in deformation 
quantization in relation with M-theory and string theory. The main hope 
seems to be that deformation quantization may  give a good framework for 
noncommutative field theory. We refer the reader to the 
thorough review by Douglas and Nekrasov~\cite{DoNe} for more details.
\vspace*{3mm}

\noindent\textbf{4.2.2 \ Path integral formula of Cattaneo and Felder.}
Kontsevich~\cite{Konsi97b} suggested that his formula comes from a perturbation
theory for a related bidimensional topological field theory.
The physical origin of the weights and graphs has been elucidated by Cattaneo 
and Felder~\cite{CaFe1}. They have constructed a topological 
field theory on a disc whose perturbation series makes Kontsevich graphs and 
weights appear explicitly after a finite renormalization.
The field model in question is the so-called Poisson sigma model~\cite{Ike};
it is made up of the following ingredients:

\noindent a) the closed unit disc $\Sigma$  in $\mathbb{R}^2$,

\noindent b)  a $d$-dimensional Poisson real manifold $(M,\alpha)$,

\noindent c)  a map $C^1$ $X\colon\Sigma\rightarrow M$,

\noindent d) a 1-form $\eta$ on $\Sigma$ taking values in the space 
of 1-forms on $M$, i.e, $\eta\in\Gamma(\Sigma, T^*\Sigma\otimes T^*M)$.

\noindent In physics, $\Sigma$ is called the worldsheet, 
and $M$ the target space.

The action functional for the Poisson sigma model is given by
the integration of the sum of two 2-forms on $\Sigma$.
Let $V_1$ and $V_2$ be two vector fields on $\Sigma$, then
$\langle \eta, V_i \rangle$, $i=1,2$, is a 1-form on $M$. 
We can define a 2-form $\beta$ on $\Sigma$ as follows:
$$
\langle\beta, V_1 \wedge V_2\rangle 
= \frac{1}{2} \big(\langle\langle \eta, V_1\rangle,DX.V_2\rangle
- \langle\langle \eta, V_2\rangle,DX.V_1\rangle \big),
$$
where $DX\colon T\Sigma\rightarrow TM$ is the tangent map of $X$.
The 2-form $\beta$ is also denoted $\langle \eta,DX\rangle$.
A second 2-form $\gamma$ on $\Sigma$ involves the Poisson tensor $\alpha$,
and is given by:
$$
\langle\gamma, V_1 \wedge V_2\rangle 
= \frac{1}{2} \langle\alpha\circ X, 
\langle \eta, V_1\rangle \wedge \langle \eta, V_2\rangle\rangle,
$$
and we shall write 
$\gamma = \frac{1}{2} \langle\alpha\circ X,\eta\wedge \eta\rangle$.
If $(u^1,u^2)$ is a chart at $s\in\Sigma$, and $(x^1,\ldots,x^d)$
is a chart at $x=X(s)$, the local expressions for $\beta$ and $\gamma$ read:
\begin{eqnarray*}
\beta(s) &=& 
\sum_{i,\mu,\nu}\eta_{\mu i}(s) 
\frac{\partial X^i}{\partial {u^\nu}} (s)\ du^\mu\wedge du^\nu\\
\gamma(s) &=&\frac{1}{2} \sum_{i,j,\mu,\nu}\alpha^{ij}(x) 
\eta_{\mu i}(s)\eta_{\nu j}(s)\ du^\mu\wedge du^\nu
\end{eqnarray*}

Finally, the action integral is defined by:
\begin{equation}\label{action}
S[X,\eta] =\int_\Sigma \big(\langle \eta,DX\rangle 
+\frac{1}{2} \langle\alpha\circ X,\eta\wedge \eta \big\rangle).
\end{equation}
The 1-form $\eta$ is required to vanish on $T(\partial\Sigma)$ (tangent vectors
of $\Sigma$ tangential to the boundary $\partial\Sigma=S^1$).

The dynamics of $S[X,\eta]$ is one of a constrained Hamiltonian system 
governing the propagation of an open string. 
Let $\alpha^{\sharp}\colon T^*M\rightarrow TM$ be the canonical 
homomorphism of vector bundles induced by $\alpha$. The variation 
of ~(\ref{action}) under the constraint $\eta|_{\partial\Sigma}=0$
gives us the critical point $(X,\eta)$ which is a solution of 
(here  $x=X(s)$):
\begin{eqnarray}
&& DX(s) + \langle\alpha^{\sharp}(x), \eta(s) \rangle =0,\nonumber\\
&& \label{motion}\\
&& d_s\eta(s) + \frac{1}{2}
d_x \langle\alpha(x),\eta(s)\wedge \eta(s)\rangle=0,\nonumber
\end{eqnarray}
where $d_s$ (resp. $d_x$) is the de Rham differential on $\Sigma$ (resp. $M$).
Equations~(\ref{motion}) admit the trivial solution: $X$ is a constant 
map equal to $x$, and $\eta=0$.

The integral formula of Cattaneo and Felder is defined through path 
integration over $(X,\eta)$ subject to the following conditions:
Take three cyclicly counterclockwise ordered points $a$, $b$, $c$ 
on the unit circle, the boundary of $\Sigma$. $X(c) = x$ a fixed point in $M$, 
$\eta|_{\partial\Sigma}=0$; then the functional integral
$$
(f\star_\hbar g)(x) =\int \mathcal{D}(X)\mathcal{D}(\eta) f(X(a)) g(X(b))
e^{\frac{i}{\hbar}S[X,\eta]}.
$$
where $f,g$ are smooth functions on $M$, ``defines'' a differential 
star-product on $(M,\alpha)$. Of course, the meaning of the integral has to 
be taken with a grain of salt, but its heuristic power cannot be denied
(here the star-product itself is expressed as a path integral, for a
zero Hamiltonian).

The proof of this formula is given by 
the quantization of this  Poisson  sigma model, which is a delicate question. 
The nature of the constraints do not allow to use BRST-methods and quantization 
is performed in the Batalin--Vilkovisky framework. We shall only say
here that the perturbative expansion of this theory provides (essentially)
the graphs and weights of Kontsevich (the angle function appears as a 
propagator). Loop graphs (or tadpole), which are excluded in Kontsevich 
formula, do appear in this perturbative expansion, but it is shown that 
they disappear by renormalization.

We conclude this section by remarking that the methods used by 
Cattaneo and Felder proved to be very efficient. For example, they
have provided a simple way to understand the globalization of 
Kontsevich star-products in more geometrical terms~\cite{CaFeTo}.
The proof goes much along the lines of Fedosov's construction.
\newpage
\subsection{Recent metamorphoses: operadic approach, algebraic varieties}

This section will deal mostly with the metamorphoses that occurred in
deformation quantization during the past three years. About one month
after the ICM98 Congress in Berlin, Tamarkin \cite{Ta98}found a new short 
derivation (for $X=\mathbb{R}^n$) of the Formality Theorem, based on a 
very general result concerning all associative algebras: For any 
algebra $A$ (over a field of characteristic 0), its cohomological 
Hochschild complex $C^\bullet(A,A)$ and its Hochschild cohomology 
$B\equiv H^\bullet(A,A)$ (this graded space 
being considered as a complex with zero differential) are algebras over 
the same operad (up to homotopy). However abstruse this statement may
seem to a physicist, when Kontsevich saw Tamarkin's result he thought to
himself, like Commissaire Maigret in Simenon's detective stories,
``Bon sang, mais c'est bien s\^ur \ldots". In fact, he had been close 
to such a result in some works of his in the early nineties but somehow
missed the point at that time. He then devoted a lot of his time 
and efforts to understand fully Tamarkin's result and generalize it, 
which he did to a considerable extent. That is probably the main reason 
why he did not finalize in print his fundamental 1997 e-print 
\cite{Konsi97b}, nor write his ICM98 lecture on the relations between 
deformations, motives (in particular, the motivic Galois group) 
and the Grothendieck--Teichm\"uller group:
all of this became much more clear when expressed in the new language. 
Instead he pushed a lot further the various ramifications of the 
new metamorphosis in two fundamental papers \cite{Konsi99,KoS}, the latter 
with Soibelman with whom he is working on what will certainly be a very
seminal book on these subjects (and many more).  

Recently he accomplished still another leap forward \cite{Konsi01} by 
going from the $C^\infty$ situation to the algebro-geometric setting,
performing what can be called `noncommutative algebraic geometry'. 
In the previous case, operads and deformations of algebras over operads
were essential notions. In the new setting he had to use sheaves of
algebroids, a slightly more complicated structure. 
We shall not attempt here to summarize in a couple of pages the very
dense content of several long and complicated papers, rather to try and
convey a few touches of (what we understand of) the flavor that emerges 
from them, giving only very basic definitions. Both developments have 
a wide array of ramifications and implications in a variety of frontier 
mathematics. 
\vspace*{3mm}

\noindent\textbf{4.3.1 \ Operadic approach.}
To fix ideas we shall give a short definition of an operad and of related 
notions as in \cite{Konsi99}. More abstract (and general) definitions 
can be found in \cite{KoS}. The language of operads, convenient
for descriptions and constructions of various algebraic structures,
became recently quite popular in theoretical physics because of the 
emergence of many new types of algebras related with quantum field theories.

\begin{definition}
An \textsf{operad} (of vector spaces) consists of the following:

\noindent O$_1$) a collection of vector spaces $P(n)$, $n\ge 0$,

\noindent O$_2$) an action of the symmetric group $S_n$ on $P(n)$ for every $n$,

\noindent O$_3$) an identity element $\mathrm{id}_P\in P(1)$,

\noindent O$_4$) compositions $m_{(n_1,\dots, n_k)}$:
\begin{equation}
P(k)\otimes \bigl(P(n_1)\otimes P(n_2)\otimes\dots\otimes 
P(n_k)\bigr) \longrightarrow P(n_1+\dots+n_k)
\end{equation} 
for every $k\ge 0$ and $n_1,\dots,n_k\ge 0$
satisfying a natural list of axioms. 

If we replace vector spaces by topological spaces, $\otimes$ by $\times$,
and require continuity of all maps, we get a \textsf{topological operad}.

An \textsf{algebra over an operad} $P$ consist of a vector space $A$ and a 
collection of polylinear maps $f_n:P(n)\otimes A^{\otimes n}\longrightarrow A$
for all $n\ge 0$ satisfying the following list of axioms:
 
\noindent  A$_1$) for any $n\ge 0$  the map  $f_n$ is $S_n$-equivariant,
  
\noindent   A$_2$) for any $a\in A$ we have $f_1(\mathrm{id}_P\otimes a)=a$,
  
\noindent   A$_3$) all compositions in $P$ map to compositions of 
polylinear operations on $A$. 

The \textsf{little discs operad} $C_d$ is a topological operad 
such that $C_d(0)=\emptyset$, 
$C_d(1)= \,\mathrm{point}\,=\{\mathrm{id}_{C_d}\}$, and
for $n\ge 2$ the space $C_d(n)$ is the space of configurations 
of $n$ disjoint discs $(D_i)_{1\le i\le n}$ inside the standard disc $D_0$.
The composition  $C_d(k)\times C_d(n_1)\times\dots\times C_d(n_k)
\longrightarrow C_d(n_1+\dots+n_k)$ is obtained by applying the affine
transformations in $AF(\mathbb{R}^d)$ associated with  discs 
$(D_i)_{1\le i\le k}$ in the configuration in $C_d(k)$ 
(the translations and dilations transforming $D_0$ to $D_i$)
to configurations in all $C_d(n_i),\,\,i=1,\dots,k$
and putting the resulting configurations together. The action of 
the symmetric group $S_n$ on $C_d(n)$ is given by renumerations of 
indices of discs $(D_i)_{1\le i\le n}$.
\end{definition}

An associative algebra is an algebra over some operad (denoted 
$\mathsf{Assoc}$ in \cite{Konsi99}). This explains the universality of
associative algebras and suggests how to develop the deformation theory
of algebras over operads, which is the right context for the deformation 
theory of many algebraic structures. 

Now, for a topological space $X$, denote by $\mathsf{Chains}(X)$ 
the complex concentrated in $\mathbb{Z}_{\le 0}$ whose $(-k)$-th component 
for $k=0,1,\dots$ consists of the formal finite additive combinations
$\sum_{i=1}^N n_i \cdot f_i$, $n_i\in \mathbb{Z}$, $N\in \mathbb{Z}_{\ge 0}$ 
of  continuous maps $f_i:[0,1]^k\longrightarrow X$ 
(singular cubes in $X$), modulo the following relations: 
$f\circ \sigma=\mathrm{sign}(\sigma) \, f$ 
for any $\sigma \in S_k$ acting on the standard cube $[0,1]^k$ by 
permutations of coordinates, and $f'\circ pr_{k\longrightarrow (k-1)}=0$ where 
$pr_{k\longrightarrow (k-1)}:[0,1]^k\longrightarrow[0,1]^{k-1}$ is the 
projection onto the first $(k-1)$ coordinates, and 
$f':[0,1]^{k-1}\longrightarrow X$ is a continuous map.
The boundary operator on cubical chains is defined in the usual way.
A $d$-algebra is an algebra over the operad $\mathsf{Chains}(C_d)$.

If $P$ is a topological operad then the collection of complexes 
$\bigl(\mathsf{Chains}(P(n))\bigr)_{n\ge 0}$ has a natural operad
structure in the category of complexes of Abelian groups. The compositions
in $\mathsf{Chains}(P)$ are defined using the external tensor product 
of cubical chains. Passing from complexes to their cohomology we obtain an 
operad $H_\bullet(P)$ of $\mathbb{Z}$-graded Abelian groups (complexes with 
zero differential), the homology operad of $P$.

A remarkable fact is the relationship of the deformation theory
of associative algebras to the geometry of configuration spaces of 
points on surfaces. One of its incarnations is the so-called Deligne 
conjecture. A formulation of the conjecture is that there is a natural action 
of $\mathsf{Chains}(C_2)$, the chain operad of the little discs operad $C_2$,
on the Hochschild complex $C^\bullet(A,A)$ for an arbitrary associative 
algebra $A$ (we refer to \cite{Konsi99,KoS} for a more precise formulation). 
The conjecture was formulated in 1993, deemed proved in 1994 until in the 
Spring of 1998 Tamarkin found a very subtle flaw in these proofs (he had 
used the result in an earlier version of his new proof of the Formality
Theorem). New proofs can be found in \cite{KoS} and other papers 
published around that time. 
The power of those notions can be grasped if we realize that the Formality
Theorem in (4.1.2) becomes a relatively easy consequence of the formality
of $\mathsf{Chains}(C_2)$, which shows that the Hochschild complex of an 
associative algebra is a homotopic Gerstenhaber algebra. A general result
is \cite{Konsi99}:
\begin{theorem} 
The operad $\mathsf{Chains}(C_d)\otimes \mathbb{R}$ of complexes of real 
vector spaces is quasi-isomorphic to its cohomology operad endowed with the
zero differential.
\end{theorem}
\noindent Then Tamarkin's result (which implies the Formality theorem; see
\cite{TaT} for further developments) can be formulated as
\begin{theorem} 
Let $A:=\mathbb{R}[x_1,\dots,x_n]$ be the algebra of polynomials considered 
merely as an associative algebra. Then the Hochschild complex 
$\mathsf{Hoch}(A)$ is quasi-isomorphic as $2$-algebra to its cohomology 
$B:=H^\bullet(\mathsf{Hoch}(A))$, the space of polynomial polyvector fields
on $\mathbb{R}^n$, considered as a Gerstenhaber algebra, hence a $2$-algebra.
\end{theorem}

Pursuing further in the directions hinted at in \cite{Konsi99}, he could
show in \cite{KoS}, among other, that the Grothendieck--Teichm\"uller group
acts homotopically on the moduli space of structures on $2$-algebras on the
Hochschild complex. The weights $w(\Gamma)$ (\ref{wa}) that appear in 
Kontsevich's quasi-isomorphism are examples of some very special numbers, 
the periods \cite{KoSMF} of rational algebraic varieties, a countable set of 
numbers lying between rational algebraic numbers and all complex numbers.
At the end of \cite{KoS} there is a theory of singular chains, suitable
for working with manifolds with corners, where one sees appearing a notion
of piecewise algebraic spaces. 
\vspace*{3mm}

\noindent\textbf{4.3.2 \ Algebraic varieties.}
There have been a number of studies of star-products in the complex
domain, starting with e.g. \cite{BG} (even if the authors had not
realized at that time that this was in fact what they were studying).
A more recent development, with a number of interesting examples, can 
be found in \cite{Bo}. However the most systematic study appeared very 
recently with Kontsevich's paper \cite{Konsi01} devoted to peculiarities
of the deformation quantization in the algebro-geometric context. 

A direct application of the Formality Theorem to an algebraic Poisson 
manifold gives a canonical sheaf of categories deforming coherent sheaves.
The global category is very degenerate in general. For a general algebraic
Poisson manifold one gets a canonical presheaf of \textsl{algebroids} 
and eventually an Abelian category, that can be quite degenerate.  
Thus, he introduces a new notion of a \textsf{semi-formal deformation}, 
a replacement in algebraic geometry of a weakened notion of strict 
deformation quantization (versus a formal one), which is quite natural 
here because of the context. To give a flavor of the kind of 
metamorphosis undergone here by deformation quantization, we reproduce 
the precise definition of \cite{Konsi01}:
\begin{definition} 
Let $R$ be  a complete pro-Artin local ring with residue field $\mathbf{k}$. 
A \textsf{deformation} over $\mathrm{Spec}(R)$ of a $\mathbf{k}$-algebra $A$ 
is an algebra $\widehat{A}$ over $R$, topologically free as $R$-module, 
considered together with an identification of 
$\mathbf{k}$-algebras $A\simeq\widehat{A}\otimes_R \mathbf{k}$. 
A \textsf{semi-formal deformation} over $R$ of a finitely generated
associative algebra $A/\mathbf{k}$ is  an algebra $\widehat{A}_\mathrm{finite}$ 
over $R$ endowed with an identification of $\mathbf{k}$-algebras 
$A\simeq \widehat{A}_\mathrm{finite}\otimes_R \mathbf{k}$ such that
there exists an exhaustive increasing filtration 
of $\widehat{A}_\mathrm{finite}$  by finitely generated
free $R$-modules $\widehat{A}_{\le n}$, compatible
with the product,  admitting a splitting as a filtration of $R$-modules, 
and such that the Rees ring of the induced filtration of $A$ is 
finitely generated over $\mathbf{k}$. 
\end{definition}

Deformed algebras obtained by semi-formal deformations are Noetherian 
and have polynomial growth. Kontsevich gives constructions of semi-formal 
quantizations (noncommutative deformations over the ring of formal
power series) of projective and affine algebraic Poisson manifolds 
satisfying certain natural geometric conditions.
Projective symplectic manifolds (e.g. K3 surfaces and abelian varieties)
do not satisfy those conditions (in contradistinction with the $C^\infty$ 
situation where quantum tori are important examples), but projective spaces 
with quadratic Poisson brackets and Poisson--Lie groups can be semi-formally 
quantized. In other words, in that framework, semi-formal deformation 
quantization is either canonical or impossible. 
\vspace*{3mm}

\noindent\textbf{4.3.3 \ Generalized deformations, singularities.}  
We shall end this (long) review by various attempts, motivated by physical 
problems, to go beyond DrG (continuous) deformations. The first two involve
deformations that are formally different from the latter, while the last
one (mentioned only for completeness) involves going to noncontinuous 
cochains in a context which has a flavor of the preceding Section 4.3.2.
\vspace{2mm}

\noindent\textsl{4.3.3.1. \ Nambu mechanics and its ``Fock--Zariski'' 
quantization.}
We mention this aspect here mainly as an example
of generalized deformation. Details can be found in \cite{DFST}, 
so we shall just briefly indicate a few highlights.

In 1973 Nambu published some calculations which he had made a dozen years 
before: with quarks in the back of his mind he started with a 
kind of ``Hamilton equations" on ${\mathbb R}^3$ with two ``Hamiltonians"
$g$, $h$ functions of ${r}$. 
In this new mechanics the evolution of a function $f$ on ${\mathbb R}^3$ is
$\frac{df}{dt}=\frac{\partial(f,g,h)}{\partial(x,y,z)}$, a 3-bracket,
where the right-hand side is the Jacobian of the mapping
${\mathbb R}^3 \rightarrow {\mathbb R}^3$ given by $(x,y,z)\mapsto (f,g,h)$.
That expression was easily generalized to $n$ functions $f_i$,
$i=1, \ldots, n$. One introduces an $n$-tuple of functions on ${\mathbb R}^n$
with composition law given by their Jacobian, linear canonical
transformations ${\rm{SL}}(n,{\mathbb R})$ and a corresponding $(n-1)$-form
which is the analogue of the Poincar\'e--Cartan integral invariant.
The Jacobian has to be interpreted as a generalized Poisson bracket:
It is skew-symmetric with respect to the $f_i$'s, satisfies an identity
(called the fundamental identity) which is an analogue of the Jacobi 
identity and is a derivation of the algebra of smooth functions on  
${\mathbb R}^n$ (i.e., the Leibniz rule is verified in each argument).
There is a complete analogy with the Poisson bracket formulation
of Hamilton equations, including the important fact that the components of the
$(n-1)$-tuple of ``Hamiltonians" $(f_2,\ldots ,f_n)$ are constants of motion.
Shortly afterwards it was shown that Nambu mechanics could
be seen as a coming from constrained Hamiltonian mechanics; e.g. for
${\mathbb R}^3$ one starts with ${\mathbb R}^6$ and an identically vanishing
Hamiltonian, takes a pair of second class constraints to reduce it to some
${\mathbb R}^4$ and one more first-class Dirac constraint, together with time
rescaling, will give the reduction. This ``chilled" the domain for almost
20 years -- and gives a physical explanation to the fact that Nambu could
not go beyond Heisenberg quantization.

In order to quantize the Nambu bracket, a natural idea is to replace, in the
definition of the Jacobian, the pointwise product of functions by a deformed
product. For this to make sense, the deformed product should be
Abelian, so we are lead to consider commutative DrG-deformations of an
associative and commutative product. Looking first at polynomials
we are lead to the commutative part of Hochschild cohomology called 
Harrison cohomology, which is trivial.
Dealing with polynomials, a natural idea is to factorize
them and take symmetrized star-products of the irreducible factors. 
More precisely we introduce an operation $\alpha$ which maps a product 
of factors into a symmetrized tensor product (in a kind of Fock space) 
and an evaluation map $T$ which replaces tensor product by star-product. 
Associativity will be satisfied if $\alpha$ annihilates the deformation 
parameter $\lambda$ (there are still $\lambda$-dependent terms in a product
due to the last action of $T$); intuitively one can think of a deformation
parameter which is $\lambda$ times a Dirac $\gamma$ matrix. 
This fact brought us to generalized deformations, but even that
was not enough. Dealing with distributivity of the product with respect 
to addition, and with derivatives, posed difficult problems. 
In the end we took for observables Taylor developments of elements
of the algebra of the semi-group generated by irreducible polynomials
(``polynomials over polynomials", inspired by second quantization techniques)
and were then able to perform a meaningful quantization of these Nambu--Poisson
brackets (cf. \cite{DFS} for more details and subsequent
developments). We call this approach ``Fock--Zariski'' to underline the 
decomposition of polynomials into irreducible factors and the role 
of these factors, similar to one-particle states in Fock quantization.
Related cohomologies were studied \cite{Ga}.
\vspace{2mm}

\noindent\textsl{4.3.3.2. \ Generalized deformations.}
The fact that in the Fock--Zariski quantization, the deformation parameter
behaves almost as if it was nilpotent, has very recently induced 
Pinczon \cite{Pi97} and Nadaud \cite{NaDef} to generalize the Gerstenhaber 
theory to the case of a deformation parameter which \textsl{does not commute 
with the algebra}. 
For instance one can have \cite{Pi97}, for $\tilde{a}=\sum_n a_n\lambda^n$,
$a_n \in A$, a left multiplication by $\lambda$ of the form 
$\lambda \cdot \tilde{a}=\sum_n \sigma(a_n)\lambda^{n+1}$ where $\sigma$ is an
endomorphism of $A$. A similar theory
can be done in this case, with appropriate cohomologies. While that theory
does not yet reproduce the above mentioned Nambu quantization, it gives
new and interesting results. In particular \cite{Pi97}, while the Weyl
algebra $W_1$ (generated by the Heisenberg Lie algebra ${\mathfrak{h}}_1$)
is known to be DrG-rigid, it can be nontrivially deformed in
such a \textsl{supersymmetric deformation theory} to the supersymmetry
enveloping algebra ${\cal{U}}({\mathfrak{osp}}(1,2))$.  
More recently \cite{NaDef}, on the polynomial algebra ${\mathbb C}[x,y]$ 
in 2 variables, Moyal-like products of a new type were discovered; a more 
general situation was studied, where the relevant Hochschild cohomology
is still valued in the algebra but with ``twists" on both sides for the 
action of the deformation parameter on the algebra.
Though this more balanced generalization of deformations also does not
(yet) give Nambu mechanics quantization, it opens a whole new direction of
research for deformation theory. This is another example of a 
physically motivated study which goes beyond a generally accepted framework 
and opens new perspectives.
\vspace{2mm}

\noindent\textsl{4.3.3.3. \ Singularities.}  We have already mentioned
the case of manifolds with boundaries or with corners, where deformation
quantization can be extended \cite{NT95}, as was the pseudodifferential 
calculus by the Melrose $b$-calculus \cite{Me}. More general situations
are being studied. An interesting new idea, where a nontrivial Harrison
cohomology can be expected, is to try and use noncontinuous cochains.
When the manifold has singularities (like a cone, a most elementary 
example), this is a reasonable thing to do \cite{CF}.
\vspace*{3mm}

\noindent\textbf{Acknowledgement.} The authors wish to thank Gilles Halbout
for his persistent and very persuasive efforts that brought them to 
eventually write down this panorama.

\frenchspacing

\end{document}